\documentclass[12pt]{article}
\usepackage{amssymb,amsmath,latexsym,amsthm}
\usepackage{mathrsfs}
\usepackage{url}
\usepackage{enumerate}
\usepackage{ifpdf}
\ifpdf
  \usepackage[pdftex]{hyperref}
\else
    \usepackage[dvips]{hyperref}
\fi

\title{Joint convergence along different subsequences of the signed cubic variation of fractional Brownian motion}
\author{Krzysztof Burdzy\thanks{Partially supported by grant DMS-1206276 from the NSF and by grant N N201 397137 from the MNiSW, Poland.}\\
  University of Washington
  \and
  David Nualart\thanks{Partially supported by grant DMS-1208625 from the NSF.}\\
  University of Kansas
  \and
  Jason Swanson\\
  University of Central Florida}
\date{October 4, 2012}

\usepackage{color}

\begin{document}

\newtheorem{thm}{Theorem}[section]
\newtheorem{cor}[thm]{Corollary}
\newtheorem{lemma}[thm]{Lemma}
\theoremstyle{definition}
\newtheorem{rmk}[thm]{Remark}
\newtheorem{expl}[thm]{Example}

\numberwithin{equation}{section}

\def\al{\alpha}
\def\be{\beta}
\def\ga{\gamma}
\def\Ga{\Gamma}
\def\de{\delta}
\def\De{\Delta}
\def\ep{\varepsilon}
\def\eps{\varepsilon}
\def\ze{\zeta}
\def\th{\theta}
\def\ka{\kappa}
\def\la{\lambda}
\def\La{\Lambda}
\def\vpi{\varpi}
\def\si{\sigma}
\def\Si{\Sigma}
\def\ph{\varphi}
\def\om{\omega}
\def\Om{\Omega}

\def\wt{\widetilde}
\def\wh{\widehat}
\def\ol{\overline}
\def\ds{\displaystyle}

\def\nab{\nabla}
\def\pa{\partial}
\def\To{\Rightarrow}
\def\eqd{\overset{d}{=}}
\def\emp{\emptyset}

\def\pf{\noindent{\bf Proof.} }
\def\qed{\hfill $\Box$}

\providecommand{\flr}[1]{\left\lfloor{#1}\right\rfloor}
\providecommand{\ceil}[1]{\left\lceil{#1}\right\rceil}
\providecommand{\ang}[1]{\left\langle{#1}\right\rangle}

\def\bA{\mathbb{A}}
\def\bB{\mathbb{B}}
\def\bC{\mathbb{C}}
\def\bD{\mathbb{D}}
\def\bE{\mathbb{E}}
\def\bF{\mathbb{F}}
\def\bG{\mathbb{G}}
\def\bH{\mathbb{H}}
\def\bI{\mathbb{I}}
\def\bJ{\mathbb{J}}
\def\bK{\mathbb{K}}
\def\bL{\mathbb{L}}
\def\bM{\mathbb{M}}
\def\bN{\mathbb{N}}
\def\bO{\mathbb{O}}
\def\bP{\mathbb{P}}
\def\bQ{\mathbb{Q}}
\def\bR{\mathbb{R}}
\def\bS{\mathbb{S}}
\def\bT{\mathbb{T}}
\def\bU{\mathbb{U}}
\def\bV{\mathbb{V}}
\def\bW{\mathbb{W}}
\def\bX{\mathbb{X}}
\def\bY{\mathbb{Y}}
\def\bZ{\mathbb{Z}}

\def\bfA{{\bf A}}
\def\bfB{{\bf B}}
\def\bfC{{\bf C}}
\def\bfD{{\bf D}}
\def\bfE{{\bf E}}
\def\bfF{{\bf F}}
\def\bfG{{\bf G}}
\def\bfH{{\bf H}}
\def\bfI{{\bf I}}
\def\bfJ{{\bf J}}
\def\bfK{{\bf K}}
\def\bfL{{\bf L}}
\def\bfM{{\bf M}}
\def\bfN{{\bf N}}
\def\bfO{{\bf O}}
\def\bfP{{\bf P}}
\def\bfQ{{\bf Q}}
\def\bfR{{\bf R}}
\def\bfS{{\bf S}}
\def\bfT{{\bf T}}
\def\bfU{{\bf U}}
\def\bfV{{\bf V}}
\def\bfW{{\bf W}}
\def\bfX{{\bf X}}
\def\bfY{{\bf Y}}
\def\bfZ{{\bf Z}}

\def\cA{\mathcal{A}}
\def\cB{\mathcal{B}}
\def\cC{\mathcal{C}}
\def\cD{\mathcal{D}}
\def\cE{\mathcal{E}}
\def\cF{\mathcal{F}}
\def\cG{\mathcal{G}}
\def\cH{\mathcal{H}}
\def\cI{\mathcal{I}}
\def\cJ{\mathcal{J}}
\def\cK{\mathcal{K}}
\def\cL{\mathcal{L}}
\def\cM{\mathcal{M}}
\def\cN{\mathcal{N}}
\def\cO{\mathcal{O}}
\def\cP{\mathcal{P}}
\def\cQ{\mathcal{Q}}
\def\cR{\mathcal{R}}
\def\cS{\mathcal{S}}
\def\cT{\mathcal{T}}
\def\cU{\mathcal{U}}
\def\cV{\mathcal{V}}
\def\cW{\mathcal{W}}
\def\cX{\mathcal{X}}
\def\cY{\mathcal{Y}}
\def\cZ{\mathcal{Z}}

\def\sA{\mathscr{A}}
\def\sB{\mathscr{B}}
\def\sC{\mathscr{C}}
\def\sD{\mathscr{D}}
\def\sE{\mathscr{E}}
\def\sF{\mathscr{F}}
\def\sG{\mathscr{G}}
\def\sH{\mathscr{H}}
\def\sI{\mathscr{I}}
\def\sJ{\mathscr{J}}
\def\sK{\mathscr{K}}
\def\sL{\mathscr{L}}
\def\sM{\mathscr{M}}
\def\sN{\mathscr{N}}
\def\sO{\mathscr{O}}
\def\sP{\mathscr{P}}
\def\sQ{\mathscr{Q}}
\def\sR{\mathscr{R}}
\def\sS{\mathscr{S}}
\def\sT{\mathscr{T}}
\def\sU{\mathscr{U}}
\def\sV{\mathscr{V}}
\def\sW{\mathscr{W}}
\def\sX{\mathscr{X}}
\def\sY{\mathscr{Y}}
\def\sZ{\mathscr{Z}}

\def\fA{\mathfrak{A}}
\def\fB{\mathfrak{B}}
\def\fC{\mathfrak{C}}
\def\fD{\mathfrak{D}}
\def\fE{\mathfrak{E}}
\def\fF{\mathfrak{F}}
\def\fG{\mathfrak{G}}
\def\fH{\mathfrak{H}}
\def\fI{\mathfrak{I}}
\def\fJ{\mathfrak{J}}
\def\fK{\mathfrak{K}}
\def\fL{\mathfrak{L}}
\def\fM{\mathfrak{M}}
\def\fN{\mathfrak{N}}
\def\fO{\mathfrak{O}}
\def\fP{\mathfrak{P}}
\def\fQ{\mathfrak{Q}}
\def\fR{\mathfrak{R}}
\def\fS{\mathfrak{S}}
\def\fT{\mathfrak{T}}
\def\fU{\mathfrak{U}}
\def\fV{\mathfrak{V}}
\def\fW{\mathfrak{W}}
\def\fX{\mathfrak{X}}
\def\fY{\mathfrak{Y}}
\def\fZ{\mathfrak{Z}}

\maketitle

\begin{abstract}

The purpose of this paper is to study the convergence in distribution of two subsequences of the signed cubic variation of the fractional Brownian motion with Hurst parameter $H=1/6$.  We prove that, under some conditions on both subsequences,  the limit is a two-dimensional Brownian motion whose components may be correlated and we find explicit formulae for its covariance function.

\bigskip
\noindent{\bf AMS subject classifications:} Primary 60G22;
secondary 60F17.

\noindent{\bf Keywords and phrases:} Fractional Brownian motion, cubic variation, convergence in law.

\end{abstract}

\section{Introduction}

Suppose that $B=\{B(t), t\ge 0\}$ is a fractional Brownian motion with Hurst parameter $H=\frac 1{2k}$, where $k$ is an odd number.  It has been proved by Nualart and Ortiz-Latorre in \cite{NualartOrtiz2008} that sequence of sums
\[
W_n(t)=\sum_{j=1}^{\flr{nt}} \left( B(j/n)-B((j-1)/n) \right) ^{k}
\]
converges in law to a  Brownian motion $W=\{W(t), t\ge 0 \}$, with variance $\sigma^2_kt$,  independent of $B$. The Brownian motion $W$ is called the $k$-signed variation of $B$. In the particular case $k=3$, the variance, denoted by $\kappa^2t$, is given in formula  (\ref{kap_def}) below. A detailed analysis of the  signed cubic variation of $B$ has been recently developed by Swanson in \cite{Swanson2011a}, considering this variation as a class of sequences of processes. 

In  the present paper, we take $H = 1/6$ and consider the case of the signed cubic variation. We are interested in the convergence in distribution of the sequence of two-dimensional processes $\{W_{a_n}(t), W_{b_n}(t)\}$, where $\{a_n\}$ and $\{b_n\}$ are two strictly increasing sequences of natural numbers converging to infinity. Under some conditions, the limit of this sequence    is a two-dimensional Gaussian process $X^\rho$, independent of $B$, whose components are Brownian motions with variance $\kappa^2t$, and with covariance $\int_0 ^t \rho(s) ds$ for some function $\rho$. The proof of this result is based on Theorem \ref{NuPe} (see Section 2.3 below), which implies that for a sequence of vectors whose components belong to a fixed Wiener chaos and each component converges in law to a Gaussian distribution, the convergence  to a multidimensional Gaussian distribution follows from the convergence of the covariance matrix. This theorem can be found in the recent monograph by Nourdin and Peccati \cite{NuPe}  (see Theorem 6.2.3) devoted to the normal approximation using  Malliavin calculus combined with Stein's method. Theorem \ref{NuPe}  has been first proved by Peccati and Tudor in   \cite{PeccatiTudor2005}, by means of stochastic calculus techniques, and Nualart and Ortiz-Latorre provide in \cite{NualartOrtiz2008}  an alternative proof based on Malliavin calculus and on the use of characteristic functions. 

The covariance function $\rho$  depends on the asymptotic behavior of the sequences $\{a_n\}$ and $\{b_n\}$.  Our main results are the following. We set $L_n =\frac{b_n} {a_n}$ and we assume that $L_n\rightarrow L\in [0,\infty]$.
\begin{itemize}
\item[(i)]  If $L=0$ of $L=\infty$, then  $\rho(s)=0$ for all $s$, and the components of $X^\rho$ are independent Brownian motions.
\item[(ii)] Suppose that $L_n=L\in(0,\infty)$ for all but finitely many $n$. Then, $L$ is a rational number, and $\rho(s)$ is a constant which depends on $L$.
\item[(iii)] If $L_n \not=L \in(0,\infty)$ for all but finitely many $n$ and  the greatest common divisor of $a_n$ and $b_n$ converges to infinity, then, again $\rho(s)$ is a constant which depends on $L$.
\item[(iv)] If $L\in(0,\infty)$ and  there exists $k\in \mathbb{N}$ such that $b_n-a_n =k$ mod $a_n$ for all $n$, then $\rho(s)$ is not constant, and depends on $L$ and $k$.
\end{itemize}
In the cases (ii)-(iv), an explicit value of $\rho(s)$ is given.

Our article is inspired by the relationship between higher (signed) variations of fractional Brownian motions and ``change of variable'' formulas for stochastic integrals with respect to these processes (see \cite{BS,NR}). These results imply that approximations to variations of fractional Brownian motion have a direct relationship with numerical stochastic integration with respect to these processes. We hope that our study will shed light on the convergence and stability of numerical approximations to stochastic integrals, and perhaps will be relevant outside the narrow context of the present article. Additionally, we find the diversity of results presented in (i)-(iv) interesting from the purely intellectual point of view, irrespective of their potential applications.

The paper is organized as follows. In Section 2 we introduce some preliminary material that will be used in the paper. We present in this section some estimates for the covariance between two increments of the fractional Brownian motion, and we study the properties of a  function $f_L(x)$, fundamental for our paper. Section 3 contains the main results and proofs, and in Section 4 we discuss some concrete examples.

\section{Preliminaries}\label{S:prelim}

If $x\in\bR$, then $\flr{x}$ denotes the greatest integer less than or equal to $x$, and $\ceil{x}$ denotes the least integer greater than or equal to $x$. Note that $\flr{x} \le x < \flr{x} + 1$, $\ceil{x} - 1 < x \le \ceil{x}$, and $\ceil{x}=\flr{x}+1_{\bZ^c}(x)$, for all $x\in\bR$. Also note that for all $n\in\bZ$ and all $x\in\bR$, we have $x<n$ if and only if $\flr{x}<n$, and $n<x$ if and only if $n <\ceil{x}$.

The Skorohod space of c\`adl\`ag functions from $[0,\infty)$ to $\bR^d$ will be denoted by $D_{\bR^d}[0,\infty)$, and convergence in law will be denoted by the symbol $\To$.

Let $B=B^{1/6}$ be a two-sided fractional Brownian motion with Hurst parameter $H=1/6$. That is, $\{B(t):t\in\bR\}$ is a centered Gaussian process with covariance function
  \[
  R(s,t) = E[B(s)B(t)]
    = \frac12(|t|^{1/3} + |s|^{1/3} - |t - s|^{1/3}),
  \]
for $s,t\in\bR$.

Let $n\in\bN$, $t_j=t_{j,n}=j/n$ and $\De B_j=\De B_{j,n} = B(t_j)-B(t_{j-1})$. If $k\in\bN$, then we shall denote $(\De B_j)^k$ by $\De B_j^k$. Let $W_n(t) = \sum_{j=1}^{\flr{nt}}\De B_j^3$. The signed cubic variation of $B$ is defined in \cite{Swanson2011a} as a class of sequences of processes, each of which is equivalent, in a certain sense, to the sequence $\{W_n\}$. The relevant fact for our present purposes is that the sequence $\{W_n\}$ converges in law to a Brownian motion independent of $B$. This was proven in \cite{NualartOrtiz2008}, and the statement of the theorem is the following.

\begin{thm}\label{T:cub_var}
As $n\to\infty$, $(B,W_n)\To(B,\ka W)$ in $D_{\bR^2}[0,\infty)$, where
  \begin{equation}\label{kap_def}
  \ka^2 = \frac34\sum_{m\in\bZ}
    (|m + 1|^{1/3} + |m - 1|^{1/3} - 2|m|^{1/3})^3,
  \end{equation}
and $W$ is a standard Brownian motion, independent of $B$.
\end{thm}

Since we are interested in the joint convergence of subsequences of $\{W_n\}$, we will be primarily concerned with the covariance of increments of this process, which can be expressed in terms of the covariance of increments of $B$. For this reason, let us define
  \begin{equation}\label{Phi_def}
  \begin{split}
  \Phi(s,t,u,v) &= 2E[(B(t) - B(s))(B(v) - B(u))]\\
  &= 2(R(t,v) - R(t,u) - R(s,v) + R(s,u))\\
  &= t^{1/3} + v^{1/3} - |t - v|^{1/3}
    - t^{1/3} - u^{1/3} + |t - u|^{1/3}\\
  &\quad - s^{1/3} - v^{1/3} + |s - v|^{1/3}
    + s^{1/3} + u^{1/3} - |s - u|^{1/3}\\
  &= |t - u|^{1/3} + |s - v|^{1/3}
    - |s - u|^{1/3} - |t - v|^{1/3},
  \end{split}
  \end{equation}
for $s,t,u,v\in\bR$. Note that
  \begin{align}
  \Phi(s,t,u,v) &= \Phi(u,v,s,t),\label{Phi_symmetry}\\
  \Phi(s,t,u,v) &= \Phi(t,t+v-u,v,v+t-s),\label{Phi_swap}\\
  \Phi(s+c,t+c,u+c,v+c) &= \Phi(s,t,u,v),\label{Phi_translation}\\
  \Phi(cs,ct,cu,cv) &= |c|^{1/3}\Phi(s,t,u,v),\label{Phi_scaling}
  \end{align}
for all $c\ge0$.

\subsection{Estimates for the function $\Phi$}\label{S:Phi}

As a first, coarse estimate of $\Phi$, note that if $x,y\in\bR$, then
  \begin{equation}
  ||x|^{1/3} - |y|^{1/3}| \le ||x| - |y||^{1/3} \le |x - y|^{1/3}.
    \label{algebra}
  \end{equation}
Thus,
  \[
  |\Phi(s,t,u,v)| \le ||t - u|^{1/3} - |s - u|^{1/3}|
    + ||s - v|^{1/3} - |t - v|^{1/3}|
    \le 2|t - s|^{1/3}.
  \]
By \eqref{Phi_symmetry},
  \[
  |\Phi(s,t,u,v)|   \le 2|v - u|^{1/3},
  \]
and it follows that
  \begin{equation}\label{basic_est}
  |\Phi(s,t,u,v)| \le 2(|t - s| \wedge |v - u|)^{1/3},
  \end{equation}
for all $s,t,u,v\in\bR$.

When more refined estimates are needed, we will rely on the following integral representations of $\Phi$. If $u<v<s<t$, then
  \begin{align}
  \Phi(s,t,u,v) &= (t - u)^{1/3} + (s - v)^{1/3}
    - (s - u)^{1/3} - (t - v)^{1/3}\notag\\
  &= \frac13\int_s^t (y - u)^{-2/3}\,dy
    - \frac13\int_s^t (y - v)^{-2/3}\,dy\notag\\
  &= -\frac29\int_s^t\int_u^v (y - x)^{-5/3}\,dx\,dy\notag\\
  &= -\frac29\int_0^{t-s}\int_0^{v-u}
    (s - v + x + y)^{-5/3}\,dx\,dy < 0.\label{int_rep1}
  \end{align}
Also, if $u<s<t<v$, then
  \begin{align}
  \Phi(s,t,u,v) &= (t - u)^{1/3} - (s - v)^{1/3}
    - (s - u)^{1/3} + (t - v)^{1/3}\notag\\
  &= \frac13\int_s^t (y - u)^{-2/3}\,dy
    + \frac13\int_s^t (y - v)^{-2/3}\,dy\notag\\
  &= \frac13\int_s^t
    ((v - y)^{-2/3} + (y - u)^{-2/3})\,dy > 0.\label{int_rep3}
  \end{align}
We will use these integral representations to generate several different estimates in Lemma \ref{L:fine_est} below.

\begin{lemma}\label{L:fine_est}
If $u<v<s<t$, then
  \begin{align}
  |\Phi(s,t,u,v)| &\le \frac29\,(t - s)(v - u)(s - v)^{-5/3},
    \label{fine_est1}\\
  |\Phi(s,t,u,v)| &\le (t - s)^{1/4}(v - u)^{11/12}(s - v)^{-5/6},
    \label{fine_est2}\\
  |\Phi(s,t,u,v)| &\le (t - s)^{11/12}(v - u)^{1/4}(s - v)^{-5/6}.
    \label{fine_est3}
  \end{align}
If $u<s<t<v$, then
  \begin{equation}
  |\Phi(s,t,u,v)| \le \frac13\,(t - s)((v - t)^{-2/3} + (s - u)^{-2/3}).
    \label{fine_est7}
  \end{equation}
\end{lemma}

\pf Suppose $u<v<s<t$. Inequality \eqref{fine_est1} follows directly from \eqref{int_rep1}. By \eqref{int_rep1} and Lemma \ref{L:inequ} (see the Appendix),
  \begin{align*}
  |\Phi(s,t,u,v)| &\le \frac29(s - v)^{-5/6}
    \int_0^{v-u} x^{-1/12}\,dx\int_0^{t-s} y^{-3/4}\,dy\\
  &= \frac29(s - v)^{-5/6} \cdot \frac{12}{11}(v - u)^{11/12}
    \cdot 4(t - s)^{1/4},
  \end{align*}
and this proves \eqref{fine_est2}. Similarly,
  \begin{align*}
  |\Phi(s,t,u,v)| &\le \frac29(s - v)^{-5/6}
    \int_0^{v-u} x^{-3/4}\,dx\int_0^{t-s} y^{-11/12}\,dy\\
  &= \frac29(s - v)^{-5/6} \cdot 4(v - u)^{1/4}
    \cdot \frac{12}{11}(t - s)^{11/12},
  \end{align*}
proving \eqref{fine_est3}. Finally, \eqref{fine_est7} follows directly from \eqref{int_rep3}. \qed

\medskip

 Let $a=\{a_n\}_{n=1}^\infty$ and $b=\{b_n\}_{n=1}^\infty$ be strictly increasing sequences in $\bN$, and let $L_n=b_n/a_n$. We define
  \begin{equation} \label{Phi}
  \Phi_n^{a,b}(j,k) = \Phi\left({
    \frac{j-1}{a_n}, \frac j{a_n},
    \frac{k-1}{b_n}, \frac k{b_n}
    }\right)
    = E[\De B_{j,a_n}\De B_{k,b_n}],
  \end{equation}
for $j,k\in\bZ$. When $a$ and $b$ are understood, we will simply write $\Phi_n$ instead of $\Phi_n^{a,b}$. By \eqref{Phi_symmetry}, we have $\Phi_n^{a,b}(j,k)=\Phi_n^{b,a}(k,j)$. Note that by \eqref{basic_est},
  \begin{equation} \label{eq2}
  |\Phi_n^{a,b}(j,k)|^3 \le 8(a_n^{-1} \wedge b_n^{-1}),
  \end{equation}
for any $a,b,n,j,k$. Applying Lemma \ref{L:fine_est} gives us Lemma \ref{L:Phi_fine}.

\begin{lemma}\label{L:Phi_fine}
If  $\frac j{a_n} < \frac {k-1}{b_n}$, then
  \begin{align}
  |\Phi_n^{a,b}(j,k)| &\le \frac29\, a_n^{-1}b_n^{-1}\left(\frac{k-1}{b_n} -\frac j {a_n} \right)^{-5/3},
    \label{Phi_fine4}\\
  |\Phi_n^{a,b}(j,k)| &\le a_n^{-1/4}b_n^{-11/12}\left(\frac{k-1}{b_n} -\frac j {a_n} \right)^{-5/6},
    \label{Phi_fine5}\\
  |\Phi_n^{a,b}(j,k)| &\le a_n^{-11/12}b_n^{-1/4}\left(\frac{k-1}{b_n} -\frac j {a_n} \right)^{-5/6}.
    \label{Phi_fine6}
  \end{align}
If $  \frac{k - 1}{b_n} < \frac{j - 1}{a_n}$ and $    \frac j{a_n} < \frac k{b_n}$,
then
  \begin{equation}\label{Phi_fine7}
  |\Phi_n^{a,b}(j,k)| \le \frac13\,a_n^{-1}
    \left({\left({\frac k{b_n} - \frac j{a_n}}\right)^{-2/3}
    + \left({\frac{j - 1}{a_n} - \frac{k - 1}{b_n}}\right)^{-2/3}
    }\right).
  \end{equation}
\end{lemma}

\subsection{The function $f_L$}\label{S:fmL}

An important function in our analysis is constructed as follows. If $m\in\bZ$ and $L\in(0,\infty)$, define $f_{m,L}\in C[0,1]$ by
  \begin{equation}\label{f_def}
  \begin{split}
  f_{m,L}(x) &= 8(E[(B(x + 1) - B(x))(B(m + L) - B(m))])^3\\
  &= \Phi(x, x + 1, m, m + L)^3\\
  &= (|x - m + 1|^{1/3} + |x - m - L|^{1/3}
    - |x - m|^{1/3} - |x - m + 1 - L|^{1/3})^3.
  \end{split}
  \end{equation}
Although $f_{m,L}(x)$ is defined only for $x\in[0,1]$, the above formula for $\Phi(x, x + 1, m, m + L)^3$ can be extended to all $x$ using \eqref{Phi_translation}. We have
  \begin{equation}\label{Phi-f_relation}
  \Phi(x, x + 1, m, m + L)^3
    = f_{m - \flr{x}, L}(x - \flr{x}),
  \end{equation}
for any $m\in\bZ$, $L\in(0,\infty)$, and $x\in\bR$. Note that by \eqref{basic_est},
  \begin{equation}\label{f_basic}
  \|f_{m,L}\|_\infty \le 8,
  \end{equation}
for any $m\in\bZ$ and $L\in(0,\infty)$. Also, by \eqref{Phi}, \eqref{Phi_scaling}, \eqref{Phi_swap}, and \eqref{Phi-f_relation},
  \begin{multline}\label{Phinfrel}
  \Phi_n^{a,b}(j,k)^3
    = \frac1{b_n}\Phi((j - 1)L_n, jL_n, k - 1, k)^3\\
    = \frac1{b_n}\Phi(jL_n, jL_n + 1, k, k + L_n)^3
    = \frac1{b_n}f_{k - \flr{jL_n},L_n}(jL_n - \flr{jL_n}).
  \end{multline}

\begin{lemma}\label{L:fkLconv}
The series $\sum_{m\in\bZ}f_{m,L}$ is absolutely convergent in $C[0,1]$ with the uniform norm.
\end{lemma}

\pf Fix $L\in(0,\infty)$. Let $m\in\bZ$ with $m<-L$. Then for any $x\in[0,1]$, we have $m < m+L < x<x+1$. Hence, by \eqref{fine_est3},
  \begin{multline*}
  |f_{m,L}(x)| = |\Phi(x,x+1,m,m+L)|^3\\
    \le L^{3/4}(x - m - L)^{-5/2}
    \le L^{3/4}(- m - L)^{-5/2}
    = L^{3/4}|m + L|^{-5/2}.
  \end{multline*}
Thus, $\|f_{m,L}\|_\infty\le L^{3/4}|m + L
 |^{-5/2}$.

Next, let $m\in\bZ$ with $m>2$. Then for any $x\in[0,1]$, we have $x<x+1<m<m+L$. Hence, by \eqref{Phi_symmetry} and \eqref{fine_est2},
  \[
  |f_{m,L}(x)| = |\Phi(x,x+1,m,m+L)|^3
    \le L^{3/4}(m - x - 1)^{-5/2}
    \le L^{3/4}|m - 2|^{-5/2}.
  \]
Thus, $\|f_{m,L}\|_\infty\le L^{3/4}|m - 2|^{-5/2}$.

Therefore, using \eqref{f_basic},
  \[
  \sum_{m\in\bZ}\|f_{m,L}\|_\infty
    \le L^{3/4}\sum_{m=-\infty}^{\ceil{-L}-1}|m + L|^{-5/2}
    + 8(3 - \ceil{-L})
    + L^{3/4}\sum_{m=3}^\infty |m - 2|^{-5/2} < \infty,
  \]
which shows that the series is absolutely convergent. \qed

\medskip

By Lemma \ref{L:fkLconv}, we may define $f_L=\sum_{m\in\bZ}f_{m,L} \in C[0,1]$. Let us also define $\wh f_{m,L}:\bR \to\bR$ by $\wh f_{m,L}(x)=f_{m,L}(x-\flr{x})$ and $\wh f_L:\bR \to\bR$ by $\wh f_L(x) =f_L(x-\flr{x})$. By Lemma \ref{L:fkLconv}, $\wh f_L=\sum_{m\in\bZ}\wh f_{m,L}$, and this series is absolutely and uniformly convergent on all of $\bR$.

In Lemma \ref{L:f_misc}, we catalog several properties of these functions that we will need later.

\begin{lemma}\label{L:f_misc}
The following relations hold:
  \begin{enumerate}[(i)]
  \item $|f_{m,L}(x) - f_{m,L'}(x)| \le 24|L - L'|^{1/3}$ for all $m$, $L$, $L'$ and $x$,
  \item if $L\in\bN$, then $f_L(x)=f_L(1-x)$ for all $x$,
  \item $f_1(1/2) < 0.1$, and
  \item $f_1(0) > 6.6$.
  \end{enumerate}
\end{lemma}

\pf Let $m\in\bZ$, $L,L'\in(0,\infty)$, and $x\in[0,1]$. By \eqref{f_basic}, we have
  \begin{multline*}
  \left|f_{m,L}(x) - f_{m,L'}(x)\right|\\
    = \left|f_{m,L}(x)^{1/3} - f_{m,L'}(x)^{1/3}\right|
   \cdot \left|f_{m,L}(x)^{2/3} + f_{m,L}(x)^{1/3}f_{m,L'}(x)^{1/3}
    + f_{m,L'}(x)^{2/3}\right|\\
    \le 12\left|f_{m,L}(x)^{1/3} - f_{m,L'}(x)^{1/3}\right|.
  \end{multline*}
Also, by \eqref{f_def} and \eqref{algebra},
  \begin{multline*}
  |f_{m,L}(x)^{1/3} - f_{m,L'}(x)^{1/3}|\\
    = ||x - m - L'|^{1/3} - |x - m - L|^{1/3}
    - |x - m + 1 - L'|^{1/3} + |x - m + 1 - L|^{1/3}|\\
    \le 2|L - L'|^{1/3}.
  \end{multline*}
Hence,
  \[
  |f_{m,L}(x) - f_{m,L'}(x)| \le 24|L - L'|^{1/3},
  \]
and this proves (i).

For (ii), let $L\in\bN$. For each $m\in\bZ$, define $\wt m=2-L-m$. Then, for all $x\in[0,1]$,
  \begin{align*}
  f_{m,L}(1 - x) &= (|x + m - 2|^{1/3} + |x + m + L - 1|^{1/3}
    - |x + m - 1|^{1/3} - |x + m - 2 + L|^{1/3})^3\\
  &= (|x - \wt m - L|^{1/3} + |x - \wt m + 1|^{1/3}
    - |x - \wt m - L + 1|^{1/3} - |x - \wt m|^{1/3})^3\\
  &= f_{\wt m,L}(x).
  \end{align*}
Since $f_L=\sum_{m\in\bZ}f_{m,L}$ and $m\mapsto\wt m$ is a bijection from $\bZ$ to $\bZ$, this proves (ii).

By \eqref{f_def}, \eqref{int_rep1}, and \eqref{fine_est1}, if $m <x-1$, then $f_{m,1}(x)<0$ and $|f_{m,1}(x)| \le \frac29(x-m-1)^{-5}$. Similarly, using \eqref{Phi_symmetry}, \eqref{int_rep1}, and \eqref{fine_est1}, if  $m>x+1$, then $f_{m,1}(x) <0$ and $|f_{m,1}(x)| \le \frac29(m-x-1)^{-5}$. It follows that
  \[
  f_1(1/2) = \sum_{m\in\bZ}f_{m,1}(1/2)
    < f_{0,1}(1/2) + f_{1,1}(1/2)
    = (3^{1/3} - 1)^3.
  \]
Since $3<24389/8000=(29/20)^3$, this gives $f_1(1/2) < (9/20)^3 = 729/8000<0.1$, proving (iii).

It also follows that
  \begin{align*}
  f_1(0) &= \sum_{m=-1}^1 f_{m,1}(0)
    - \sum_{\substack{m\in\bZ\\|m|\ge2}} |f_{m,1}(0)|\\
  &= 8 - 2(2 - 2^{1/3})^3
    - \sum_{\substack{m\in\bZ\\|m|\ge2}} |f_{m,1}(0)|\\
  &\ge 8 - 2(2 - 2^{1/3})^3
    - \frac29\sum_{\substack{m\in\bZ\\|m|\ge2}} ||m| - 1|^{-5}\\
  &= 8 - 2(2 - 2^{1/3})^3
    - \frac49\sum_{m=1}^\infty m^{-5}.
  \end{align*}
By Lemma \ref{L:numeric},
  \[
  \sum_{m=1}^\infty m^{-5} = 1 + \sum_{m=2}^\infty m^{-5}
    \le \frac54.
  \]
Thus, since $2>125/64=(5/4)^3$, we have
  \[
  f_1(0) > \frac{67}9 - 2\left({2 - \frac54}\right)^3
    = \frac{1901}{288} > \frac{1900.8}{288}
    = 6.6,
  \]
and this proves (iv). \qed

\subsection{Convergence in law of random vectors in a fixed Wiener chaos}

We denote by $\mathcal{H}(B)$ the closed linear subspace of $L^2(\Omega)$ generated by the family of random variables $\{B(t), t\ge 0\}$. For each integer $q\ge 1$, we denote by $\mathcal{H}_q$ the $q$-Wiener chaos defined as the subspace of $L^2(\Omega)$ spanned by the random variables $\{h_q(F), F\in \mathcal{H}(B), E(F^2)=1\}$, where
  \[
  h_q(x) = (-1)^q e^{x^2/2}\frac {d^q}{dx^q}(e^{-x^2/2})
  \]
 is the $q$th Hermite polynomial. Notice that $\mathcal{H}_1 =\mathcal{H}(B)$.

We finish this section with a result on the convergence of vectors whose components belong to a fixed Wiener chaos (see Theorem  6.2.3 in \cite{NuPe}).  

\begin{thm} \label{NuPe}
Let $d\ge 2$ and $q_1, \dots, q_d \ge 1$ be some fixed integers. Consider the sequence of vectors 
$F_n=(F_{1,n}, \dots, F_{d,n})$, where for each $i=1,\dots, d$, each component $F_{i,n}$ belongs to the Wiener chaos $\mathcal{H}_{q_i}$.  Suppose that
\[
\lim_{n\rightarrow \infty} E[F_{i,n} F_{j,n}]=C(i,j), \quad 1\le i,j \le d,
\]
where $C$ is a symmetric non-negative definite matrix. Then, the following two conditions are equivalent:
\begin{itemize}
\item[(i)]  $F_n$ converges in law to a $d$-dimensional Gaussian distribution $N(0,C)$.
\item[(ii)] For each $1\le i \le d$, $F_{i,n}$ converges in law to $N(0,C(i,i))$.
\end{itemize}
\end{thm}

\section{Main results and proofs}

Recall from Section \ref{S:prelim} that $W_n(t) = \sum_{j=1}^{\flr{nt}} \De B_j^3$, and that $(B,W_n)\To(B,\ka W)$, where $W$ is a Brownian motion. We wish to investigate the joint convergence in law of $(B, W_{a_n}, W_{b_n})$, where $\{W_{a_n}\}$ and $\{W_{b_n}\}$ are two different subsequences of $\{W_n\}$. Our first theorem, Theorem \ref{T:main1}, reduces this to an investigation of the asymptotic covariance.

\begin{thm}\label{T:main1}
Let $\{a_n\}_{n=1}^\infty$ and $\{b_n\}_{n=1}^\infty$ be strictly increasing sequences in $\bN$. Let $\rho\in C[0,\infty)$, and suppose that
  \begin{equation}\label{main1}
  \lim_{n\to\infty}E[W_{a_n}(s)W_{b_n}(t)]
    = \int_0^{s\wedge t}\rho(x)\,dx,
  \end{equation}
for all $0\le s,t<\infty$. Then $\|\rho\|_\infty\le\ka^2$, and we may define
  \begin{equation}\label{main1a}
  \si(t) = \ka\begin{pmatrix}
    \sqrt{1 - |\ka^{-2}\rho(t)|^2} &\ka^{-2}\rho(t)\\
    0 &1
    \end{pmatrix}.
  \end{equation}
Let $W$ be a standard, 2-dimensional Brownian motion, independent of $B$, and define
  \begin{equation}\label{main1b}
  X^\rho(t) = \int_0^t \si(s)\,dW(s).
  \end{equation}
Then $(B,W_{a_n},W_{b_n})\to(B,X^\rho)$ in law in $D_{\bR^3}[0,
\infty)$ as $n\to\infty$.
\end{thm}

\begin{rmk}
We know from \cite{NualartOrtiz2008} that $E|W_n(t)-W_n(s)|^2\to\ka^2|t-s|$. Thus, if \eqref{main1} is satisfied for some continuous $\rho$, then by H\"older's inequality,
  \[
  \int_s^t\rho(x)\,dx = \lim_{n\to\infty}
    E[(W_{a_n}(t) - W_{a_n}(s))(W_{b_n}(t) - W_{b_n}(s))]
    \le \ka^2(t - s),
  \]
for all $s<t$. Since $\rho$ is continuous, this implies $\|\rho\|_{\infty}\le \ka^2$, so that $\si(t)$ is well-defined by \eqref{main1a}. \qed
\end{rmk}

\begin{rmk}\label{R:main}
For any $j=1,\dots, n$, the random variable  $\Delta B_j^3$  can be expressed as
\[
\Delta B_j^3= n^{-1/2}h_3(n^{1/6} \Delta B_j) + 3 n^{-1/3} \Delta B_j,
\]
where $h_3(x)= x^3-3x$ is the third Hermite polynomial.  Define
\begin{equation} \label{wt}
\wt W_n(t)=W_n(t)-3n^{-1/3}B(\flr{nt}/n)=\sum_{j=1}^{\flr{nt}} n^{-1/2}h_3(n^{1/6}\De B_j).
\end{equation}
Then, for any $p\ge 2$ and any $t\ge 0$,
\begin{equation} \label{lim}
\lim_{n\rightarrow \infty} \sup_{0\le s \le t}E[|  \wt W_n(s)-W_n(s)|^p]=0.
\end{equation}
\qed
\end{rmk}

\noindent \textbf{Proof of Theorem \ref{T:main1}.}  Taking into account (\ref{lim}), it  suffices to establish the desired limit theorem for the sequence of processes $X_n= (X^1_n, X^2_n, X^3_n) := (B, \wt W_{a_n}, \wt W_{b_n})$.
 
We know  (see, for instance, \cite{NualartOrtiz2008}) that this sequence is tight in $(D_{\bR}[0,\infty))^3$. It is well-known (see Lemma 2.2 in \cite{Nourdin2010}, for example) that since the limit processes are continuous, this implies the sequence is tight in $D_{\mathbb{R}^3} [0,\infty)$. Thus to show the convergence in law it suffices to establish the convergence in law of the finite dimensional distributions. Consider a finite set of times $0\le t_1 <t_2< \cdots <t_M$ and the $3M$-dimensional random vector $ (X_n(t_1), \dots, X_n(t_M))$. This sequence of vectors satisfies the following properties:
 \begin{enumerate}
 \item The components  $X_n^i(t_j)$, $1\le i\le 3$, $1\le j \le M$, belong to the third Wiener chaos if $i=2,3$ and to the first Wiener chaos if $i=1$. 
  \item The first component $X_n^1(t_j)=B(t_j)$ is Gaussian with a fixed law.
 On the other  hand, we know from \cite{NualartOrtiz2008} that the other two components $X_n^2(t_j)= \wt W_{a_n}(t_j)$ and  $X_n^3(t_j)= \wt W_{b_n}(t_j)$ converge in law as $n$ tends to infinity to  a Gaussian distribution with variance $\kappa^2 t_j$, which coincides with the common law of   $X^{\rho,1}(t_j)$ and $X^{\rho,2}(t_j)$.
 \end{enumerate}
 
Set  $X=(B, X^{\rho,1}, X^{\rho,2})$. Then,  by Theorem \ref{NuPe}, in order to show that
  \[
  (X_n(t_1), \dots, X_n(t_M)) \Rightarrow (X(t_1), \dots, X(t_M)),
  \]
 it suffices to show that for any $i\not =k$  and for any $s,t \ge 0$, we have
\begin{equation} \label{eq1}
\lim_{n\rightarrow \infty} E[X_n^i(s)  X_n^k(t)]=E[X^i(s)X^k(t)].
\end{equation}
If $i=1$ and $k=2,3$, then  $E[X^i(s)X^k(t)]=0$ and (\ref{eq1}) has been proved in \cite{NualartOrtiz2008}. For $i=2$ and $k=3$, then, taking into account  (\ref{lim}) and using our assumption  (\ref{main1}) we obtain
\[
\lim_{n\rightarrow \infty} E[X_n^2(s)  X_n^3(t)]=
\lim_{n\rightarrow \infty} E[W_{a_n}(s)W_{b_n}(t)]=
\int_0^{s\wedge t} \rho(u) du = E[X^2(s)X^3(t)],
\]
and the proof is complete.
\qed

\medskip

Our first main result concerns the simplest of the situations we consider, where $\{W_{a_n}\}$ and $\{W_{b_n}\}$ are subsequences such that $b_n/a_n$ converges to either $0$ or $\infty$.

\begin{thm}\label{T:main2'}
Let $\{a_n\}_{n=1}^\infty$ and $\{b_n\}_{n=1}^\infty$ be strictly increasing sequences in $\bN$. Let $L_n=b_n/a_n$ and suppose that $L_n\to L\in\{0,\infty\}$. Then
    \[
    \lim_{n\to\infty}
      E[(W_{a_n}(t) - W_{a_n}(s))(W_{b_n}(t) - W_{b_n}(s))] = 0,
    \]
for all $0\le s\le t$.
\end{thm}

\pf 
 From (\ref{lim}), it suffices to show that
  \[
  \lim_{n\to\infty}
    E[(\wt W_{a_n}(t) - \wt W_{a_n}(s))
    (\wt W_{b_n}(t) - \wt W_{b_n}(s))] = 0,
  \]
where $\wt W$ has been defined in (\ref{wt}). By interchanging the roles of $\{a_n\}$ and $\{b_n\}$ if necessary, we may assume that $L=0$. Fix $0\le s\le t$ and note that
  \begin{align*}
  E[(\wt W_{a_n}(t) - \wt W_{a_n}(s))&(\wt W_{b_n}(t) - \wt W_{b_n}(s))]\\
    &= \sum_{j=\flr{a_ns}+1}^{\flr{a_nt}}\sum_{k=\flr{b_ns}+1}^{\flr{b_nt}}
    n^{-1}E[h_3(n^{1/6}\De B_{j,a_n})h_3(n^{1/6}\De B_{k,b_n})]\\
  &= \sum_{j=\flr{a_ns}+1}^{\flr{a_nt}}\sum_{k=\flr{b_ns}+1}^{\flr{b_nt}}
    n^{-1}3!(E[n^{1/6}\De B_{j,a_n}n^{1/6}\De B_{k,b_n}])^3\\
  &= 3!\sum_{j=\flr{a_ns}+1}^{\flr{a_nt}}\sum_{k=\flr{b_ns}+1}^{\flr{b_nt}}
    (E[\De B_{j,a_n}\De B_{k,b_n}])^3\\
  &= \frac34\sum_{j=\flr{a_ns}+1}^{\flr{a_nt}}\sum_{k=\flr{b_ns}+1}^{\flr{b_nt}}
    \Phi_n(j,k)^3,
  \end{align*}
where $\Phi_n(j,k)=\Phi_n^{a,b}(j,k)$ has been introduced in (\ref{Phi}). Note that in the second equality above, we have used the fact that if $X$ and $Y$ are jointly Gaussian, each with mean zero and variance one, then $E[h_q(X)h_q(Y)]=q!(E[XY])^q$. See \cite[Lemma 1.1.1]{Nualart2006}.

Define
  \[
  S_n := \frac34\sum_{j=1}^{\flr{a_nt}}\sum_{k=1}^{\flr{b_nt}}
    |\Phi_n(j,k)|^3
    \ge |E[(\wt W_{a_n}(t) - \wt W_{a_n}(s))
    (\wt W_{b_n}(t) - \wt W_{b_n}(s))]|,
  \]
so that it will suffice to show that $S_n\to0$ as $n\to\infty$.

For each fixed $k\in\{1,\ldots,\flr{b_nt}\}$, consider the following sets of indices:
  \begin{align*}
  A^k_1 &= \left\{1 \le j \le \flr{a_nt}:
    \frac{j-1}{a_n} \le \frac{k-1}{b_n}\right\},\\
  A^k_2 &= \left\{1 \le j \le \flr{a_nt}:
    \frac{k-1}{b_n} < \frac{j-1}{a_n} < \frac j{a_n} < \frac k{b_n}\right\},\\
  A^k_3 &= \left\{1 \le j \le \flr{a_nt}:
    \frac k{b_n} \le \frac j{a_n}\right\}.
  \end{align*}
It is easily verified that $\bigcup_\ell A^k_\ell=\{1,\ldots,\flr{a_nt}\}$, $A^k_1\cap A^k_2=\emptyset$, and $A^k_2\cap A^k_3=\emptyset$. Also, if $L_n<1$, then $A^k_1\cap A^k_3=\emptyset$. Thus, for $n$ sufficiently large, $\{A^k_1, A^k_2, A^k_3\}$ is a partition of $\{1,\ldots,\flr{a_nt}\}$, and we may write
  \[
  S_n = S_n^{(1)} + S_n^{(2)} + S_n^{(3)},
  \]
where
  \[
  S_n^{(i)} = \sum_{k=1}^{\flr{b_nt}}\sum_{j\in A_i}|\Phi_n(j,k)|^3.
  \]
Note that $(j-1)/a_n\le(k-1)/b_n$ if and only if $j\le\flr{(k-1)/L_n}+1$, and $j/a_n<(k-1)/b_n$ if and only if $j<\ceil{(k-1)/L_n}$. Also note that for $n$ sufficiently large, $\flr{(k-1)/L_n}+1\le\flr{a_nt}$. Thus, by \eqref{eq2} and \eqref{Phi_fine6},
  \begin{align*}
  \sum_{j\in A_1}|\Phi_n(j,k)|^3
    &= \sum_{j=1}^{\ceil{(k-1)/L_n}-3}|\Phi_n(j,k)|^3
    + \sum_{j=\ceil{(k-1)/L_n}-2}^{\flr{(k-1)/L_n}+1}|\Phi_n(j,k)|^3\\
  &\le \sum_{j=1}^{\ceil{(k-1)/L_n}-3}
    a_n^{-11/4}b_n^{-3/4}\left({
    \frac{k-1}{b_n} - \frac j{a_n}}\right)^{-5/2}
    + 32(a_n^{-1}\wedge b_n^{-1})\\
  &\le a_n^{-1/4}b_n^{-3/4}\sum_{j=1}^{\ceil{(k-1)/L_n}-3}
    \left({\frac{k-1}{L_n} - j}\right)^{-5/2}
    + 32a_n^{-1}.
  \end{align*}
Using Lemma \ref{L:numeric}, we have
  \begin{align*}
  \sum_{j\in A_1}|\Phi_n(j,k)|^3
    &\le \frac23a_n^{-1/4}b_n^{-3/4}\left({
    \frac{k-1}{L_n} - \ceil{\frac{k-1}{L_n}} + 2}\right)^{-3/2}
    + 32a_n^{-1}\\
  &\le \frac23a_n^{-1/4}b_n^{-3/4} + 32a_n^{-1},
  \end{align*}
giving
  \[
  0 \le S_n^{(1)} \le b_nt\left({
    \frac23a_n^{-1/4}b_n^{-3/4} + 32a_n^{-1}}\right)
    \le t\left({\frac23L_n^{1/4} + 32L_n}\right)
    \to 0,
  \]
as $n\to\infty$.

For $S_n^{(3)}$, note that $k/b_n\le j/a_n$ if and only if $\ceil{k/L_n}\le j$. Also, $k/b_n<(j-1)/a_n$ if and only if $j>\flr{k/L_n}+1$. Since $\Phi_n^{a,b}(j,k) =\Phi_n^{b,a}(k,j)$, we apply \eqref{Phi_fine5} with $j,k$ and $a,b$ interchanged. Using also \eqref{eq2}, we obtain
  \begin{align*}
  \sum_{j\in A_3}|\Phi_n(j,k)|^3
    &= \sum_{j=\ceil{k/L_n}}^{\flr{k/L_n}+3}|\Phi_n(j,k)|^3
    + \sum_{j=\flr{k/L_n}+4}^{\flr{a_nt}}|\Phi_n(j,k)|^3\\
  &\le 32(a_n^{-1}\wedge b_n^{-1})
    + \sum_{j=\flr{k/L_n}+4}^{\flr{a_nt}}
    a_n^{-11/4}b_n^{-3/4}\left({
    \frac{j-1}{a_n} - \frac k{b_n}}\right)^{-5/2}\\
  &\le 32a_n^{-1}
    + a_n^{-1/4}b_n^{-3/4}\sum_{j=\flr{k/L_n}+4}^{\flr{a_nt}}
    \left({j - 1 - \frac k{L_n}}\right)^{-5/2}.
  \end{align*}
By Lemma \ref{L:numeric}, we have
  \begin{align*}
  \sum_{j\in A_3}|\Phi_n(j,k)|^3
    &\le 32a_n^{-1}
    + \frac23a_n^{-1/4}b_n^{-3/4}\left({
    \flr{\frac k{L_n}} + 2 - \frac k{L_n}}\right)^{-3/2}\\
  &\le 32a_n^{-1} + \frac23a_n^{-1/4}b_n^{-3/4},
  \end{align*}
giving, as above, $S_n^{(3)}\to 0$ as $n\to\infty$.

Finally, for $S_n^{(2)}$, note that for sufficiently large $n$, we have $L_n < 1$, which implies $b_n^{-1} - a_n^{-1}>0$, giving
  \[
  \left({\frac k{b_n} - \frac j{a_n}}\right)^{-2/3}
    < \left({\frac {k-1}{b_n} - \frac {j-1}{a_n}}\right)^{-2/3}
    = \left({\frac {j-1}{a_n} - \frac {k-1}{b_n}}\right)^{-2/3}.
  \]
Hence, by \eqref{Phi_fine7},
  \begin{align*}
  \sum_{j\in A_2}|\Phi_n(j,k)|^3
    &= \sum_{j=\flr{(k-1)/L_n}+2}^{\ceil{k/L_n}-1}|\Phi_n(j,k)|^3\\
  &\le 16(a_n^{-1}\wedge b_n^{-1})
    + \frac8{27}a_n^{-3}\sum_{j=\flr{(k-1)/L_n}+4}^{\ceil{k/L_n}-1}
    \left({\frac {j-1}{a_n} - \frac {k-1}{b_n}}\right)^{-2}\\
  &\le 16a_n^{-1}
    + a_n^{-1}\sum_{j=\flr{(k-1)/L_n}+4}^{\ceil{k/L_n}-1}
    \left({j - 1 - \frac{k-1}{L_n}}\right)^{-2}.
  \end{align*}
By Lemma \ref{L:numeric}, we have
  \[
  \sum_{j\in A_2}|\Phi_n(j,k)|^3
    \le 16a_n^{-1} + a_n^{-1}\left({
    \flr{\frac{k-1}{L_n}} + 2 - \frac{k-1}{L_n}}\right)^{-1}
    \le 17a_n^{-1},
  \]
giving, as above, $S_n^{(2)}\to 0$ as $n\to\infty$. \qed

\medskip

To use Theorem \ref{T:main1}, we must verify hypothesis \eqref{main1}. Our next lemma, Lemma \ref{L:main1}, simplifies this task, allowing us to check \eqref{main1} only when $s=t$.

\begin{lemma}\label{L:main1}
Let $\{a_n\}_{n=1}^\infty$ and $\{b_n\}_{n=1}^\infty$ be strictly increasing sequences in $\bN$. Let $L_n=b_n/a_n$ and suppose that $L_n\to L\in[0,\infty]$. Then
  \[
  \lim_{n\to\infty}E[W_{a_n}(s)(W_{b_n}(t) - W_{b_n}(s))]
    = \lim_{n\to\infty}E[(W_{a_n}(t) - W_{a_n}(s))W_{b_n}(s)] = 0,
  \]
for any $0\le s<t$.
\end{lemma}

\pf By interchanging the roles of $\{a_n\}$ and $\{b_n\}$ if necessary, we may assume that $L_n\to L\in(0,\infty]$. From (\ref{lim}), it suffices to show
  \begin{equation}\label{E:main1.01}
  \lim_{n\to\infty}E[\wt W_{a_n}(s)(\wt W_{b_n}(t) - \wt W_{b_n}(s))] = 0,
  \end{equation}
and
  \begin{equation}\label{E:main1.02}
  \lim_{n\to\infty}E[(\wt W_{a_n}(t) - \wt W_{a_n}(s))\wt W_{b_n}(s)] = 0,
  \end{equation}
where $\wt W$ has been defined in (\ref{wt}). We begin by proving \eqref{E:main1.01}.

As in the proof of Theorem \ref{T:main2'}, we have
  \begin{equation}\label{main1.4}
  E[\wt W_{a_n}(s)(\wt W_{b_n}(t) - \wt W_{b_n}(s))]
    = \frac34\sum_{j=1}^{\flr{a_ns}}\sum_{k=\flr{b_ns}+1}^{\flr{b_nt}}
    \Phi_n(j,k)^3.
  \end{equation}
We claim that for all $i\ge0$,
  \begin{equation}\label{main1.3}
  \sum_{k=\flr{b_ns}+1}^{\flr{b_nt}}
    \Phi_n(\flr{a_ns}-i,k)^3 \to 0,
  \end{equation}
as $n\to\infty$. By \eqref{eq2}, it is enough to show that
  \[
  \sum_{k=\flr{b_ns}+3}^{\flr{b_nt}}
    \Phi_n(\flr{a_ns}-i,k)^3 \to 0.
  \]
For this, fix $n$ and let $j=\flr{a_ns}-i$. Note that since $\flr{x}\le x<\flr{x}+1$, we have
  \begin{equation}\label{main1.2}
  \frac j{a_n} \le \frac{a_ns - i}{a_n} \le s = \frac{b_ns}{b_n}
    < \frac{\flr{b_ns} + 1}{b_n} \le \frac{k - 2}{b_n},
  \end{equation}
for any $k\ge\flr{b_ns}+3$. Hence, by \eqref{Phi_fine5} we have
  \begin{equation}\label{main1.6}
  |\Phi_n(j,k)|  \le a_n^{-1/4}b_n^{-11/12}
    \left({\frac{k - 1}{b_n} - \frac j{a_n}}\right)^{-5/6}.
  \end{equation}
Using \eqref{main1.2}, this gives
  \begin{align*}
  \sum_{k=\flr{b_ns}+3}^{\flr{b_nt}}
    |\Phi_n(\flr{a_ns}-i,k)|^3
    &\le a_n^{-3/4}b_n^{-11/4}\sum_{k=\flr{b_ns}+3}^{\flr{b_nt}}
    \left({\frac{k - 1}{b_n}
    - \frac{\flr{b_ns} + 1}{b_n}}\right)^{-5/2}\\
  &= a_n^{-3/4}b_n^{-1/4}\sum_{k=\flr{b_ns}+3}^{\flr{b_nt}}
    (k - \flr{b_ns} - 2)^{-5/2}\\
  &\le a_n^{-3/4}b_n^{-1/4}\sum_{k=1}^{\infty} k^{-5/2}
    \to 0,
  \end{align*}
as $n\to\infty$, and this prove \eqref{main1.3}.

Now, since $b_n/a_n\to L\in(0,\infty]$, there exists an integer $\ell\ge 2$ such that $b_n/a_n\ge1/(\ell-1)$ for all $n$. We next claim that for all $i\ge0$,
  \begin{equation}\label{main1.5}
  \sum_{j=1}^{\flr{a_ns}-\ell}
    \Phi_n(j,\flr{b_ns}+i)^3 \to 0,
  \end{equation}
as $n\to\infty$. Again, fix $n$ and let $k=\flr{b_ns}+i$. Then, for all $j\le\flr{a_ns}-\ell$, we have
  \begin{equation}\label{main1.7}
  \frac j{a_n} < \frac{\flr{a_ns} - \ell + 1}{a_n}
    \le s - \frac{\ell - 1}{a_n} \le s - \frac1{b_n}
    = \frac{b_ns - 1}{b_n} < \frac{\flr{b_ns}}{b_n}
    = \frac{k - i}{b_n} \le \frac{k - 1}{b_n}.
  \end{equation}
Since $j/a_n<(k-1)/b_n$, from \eqref{Phi_fine6} we conclude
  \[
  |\Phi_n(j,k)| \le a_n^{-11/12}b_n^{-1/4}
    \left({\frac{k - 1}{b_n} - \frac j{a_n}}\right)^{-5/6}.
  \]
Using \eqref{main1.7}, this gives
  \begin{align*}
  \sum_{j=1}^{\flr{a_ns}-\ell} |\Phi_n(j,\flr{b_ns}+i)|^3
    &\le a_n^{-11/4}b_n^{-3/4}\sum_{j=1}^{\flr{a_ns}-\ell}
    \left({\frac{\flr{a_ns} - \ell + 1}{a_n}
    - \frac j{a_n}}\right)^{-5/2}\\
  &= a_n^{-1/4}b_n^{-3/4}\sum_{j=1}^{\flr{a_ns}-\ell}
    (\flr{a_ns} - \ell + 1 - j)^{-5/2}\\
  &\le a_n^{-1/4}b_n^{-3/4}\sum_{j=1}^\infty j^{-5/2}
    \to 0,
  \end{align*}
as $n\to\infty$, and this proves \eqref{main1.5}.

Finally, \eqref{E:main1.01} will be proved once we show that the double sum in \eqref{main1.4} converges to zero. Let us write
  \begin{align*}
  \sum_{j=1}^{\flr{a_ns}}
    \sum_{k=\flr{b_ns}+1}^{\flr{b_nt}} \Phi_n(j,k)^3
    &= \sum_{i=0}^{\ell-1}\sum_{k=\flr{b_ns}+1}^{\flr{b_nt}}
    \Phi_n(\flr{a_ns}-i,k)^3
    + \sum_{j=1}^{\flr{a_ns}-\ell}
    \sum_{k=\flr{b_ns}+1}^{\flr{b_nt}} \Phi_n(j,k)^3\\
  &= \sum_{i=0}^{\ell-1}\sum_{k=\flr{b_ns}+1}^{\flr{b_nt}}
    \Phi_n(\flr{a_ns}-i,k)^3
    + \sum_{i=1}^3\sum_{j=1}^{\flr{a_ns}-\ell}
    \Phi_n(j,\flr{b_ns}+i)^3\\
  &\qquad + \sum_{j=1}^{\flr{a_ns}-\ell}
    \sum_{k=\flr{b_ns}+4}^{\flr{b_nt}} \Phi_n(j,k)^3.
  \end{align*}
By \eqref{main1.3} and \eqref{main1.5}, the first two double sums above converge to zero. Hence, it will suffice to show that
  \[
  \ep_n := \sum_{j=1}^{\flr{a_ns}-\ell}
    \sum_{k=\flr{b_ns}+4}^{\flr{b_nt}} \Phi_n(j,k)^3 \to 0,
  \]
as $n\to\infty$.

As before, for all $j\le\flr{a_ns}-\ell$ and all $k\ge\flr{b_ns}+4$, we have
  \[
  \frac j{a_n} \le s - \frac\ell{a_n} < \frac{b_ns - 1}{b_n}
    \le \frac{\flr{b_ns}}{b_n} \le \frac{k - 4}{b_n},
  \]
and the estimate (\ref{Phi_fine4}) implies 
  \[
  |\Phi_n(j,k)| \le \frac29a_n^{-1}b_n^{-1}
    \left({\frac{k-1}{b_n} - \frac j{a_n}}\right)^{-5/3},
  \]
so that
  \begin{align*}
  |\ep_n| &\le a_n^{-3}b_n^{-3}\sum_{j=1}^{\flr{a_ns}-\ell}
    \sum_{k=\flr{b_ns}+4}^{\flr{b_nt}}
    \left({\frac{k-1}{b_n} - \frac j{a_n}}\right)^{-5}\\
  &= a_n^2b_n^{-3}\sum_{k=\flr{b_ns}+4}^{\flr{b_nt}}
    \sum_{j=1}^{\flr{a_ns}-\ell}
    \left({\frac{(k-1)a_n}{b_n} - j}\right)^{-5}.
  \end{align*}
To apply Lemma \ref{L:numeric}, we check that
  \[
  \flr{a_ns} - \ell + 1 < a_ns + \frac{a_n}{b_n}
    = \frac{(b_ns + 1)a_n}{b_n}
    < \frac{(\flr{b_ns} + 2)a_n}{b_n}
    < \frac{(k - 1)a_n}{b_n}.
  \]
Thus,
  \begin{align*}
  |\ep_n| &\le a_n^2b_n^{-3}\sum_{k=\flr{b_ns}+4}^{\flr{b_nt}}
    \left({\frac{(k-1)a_n}{b_n}
    - (\flr{a_ns} - \ell + 1)}\right)^{-4}\\
  &= a_n^{-2}b_n\sum_{k=\flr{b_ns}+4}^{\flr{b_nt}}\left({k - 1
    - \frac{(\flr{a_ns} - \ell + 1)b_n}{a_n}}\right)^{-4}.
  \end{align*}
To apply Lemma \ref{L:numeric} once again, we check that
  \[
  1 + \frac{(\flr{a_ns} - \ell + 1)b_n}{a_n}
    < 1 + \frac{\flr{a_ns}b_n}{a_n}
    \le 1 + b_ns < \flr{b_ns} + 3,
  \]
so that
  \begin{align*}
  |\ep_n| &\le a_n^{-2}b_n\left({\flr{b_ns} + 2
    - \frac{(\flr{a_ns} - \ell + 1)b_n}{a_n}}\right)^{-3}\\
  &= a_n^{-2}b_n^{-2}\left({\frac{\flr{b_ns} + 2}{b_n}
    - \frac{\flr{a_ns} - \ell + 1}{a_n}}\right)^{-3}.
  \end{align*}
Recalling that $\ell\ge 2$, this gives
  \[
  |\ep_n| < a_n^{-2}b_n^{-2}\left({\frac{b_ns + 1}{b_n}
    - \frac{a_ns - 1}{a_n}}\right)^{-3}
    = a_n^{-2}b_n^{-2}\left({\frac1{a_n} + \frac1{b_n}}\right)^{-3}.
  \]
Finally, by Lemma \ref{L:inequ}, we have
  \[
  |\ep_n| \le a_n^{-2}b_n^{-2}\left({\frac1{a_n}}\right)^{-3/2}
    \left({\frac1{b_n}}\right)^{-3/2}
    = a_n^{-1/2}b_n^{-1/2} \to 0,
  \]
as $n\to\infty$, and this concludes the proof of \eqref{E:main1.01}.

For \eqref{E:main1.02}, note that
  \begin{multline*}
  E[(\wt W_{a_n}(t) - \wt W_{a_n}(s))\wt W_{b_n}(s)]
    = E[\wt W_{a_n}(t)\wt W_{b_n}(t)]
    - E[\wt W_{a_n}(s)\wt W_{b_n}(s)]\\
    - E[(\wt W_{a_n}(t) - \wt W_{a_n}(s))(\wt W_{b_n}(t) - \wt W_{b_n}(s))]
    - E[\wt W_{a_n}(s)(\wt W_{b_n}(t) - \wt W_{b_n}(s))].
  \end{multline*}
By Theorem \ref{T:main2'} and Remark \ref{R:main}, the first three expectations on the right-hand side tend to zero; and by \eqref{E:main1.01}, the fourth expectation on the right-hand side tends to zero. This proves \eqref{E:main1.02} and completes the proof of the lemma. \qed

\medskip

As a consequence of  Theorem \ref{T:main1}, Lemma \ref{L:main1}, and Theorem \ref{T:main2'} we obtain the following limit result in the case $L\in \{0,\infty\}$.

\begin{cor}\label{C:main3'}
Let $\{a_n\}_{n=1}^\infty$ and $\{b_n\}_{n=1}^\infty$ be strictly increasing sequences in $\bN$. Let $L_n=b_n/a_n$ and suppose that $L_n\to L\in\{0,\infty\}$. Then $(B,W_{a_n},W_{b_n}) \to (B,\ka W^1, \ka W^2)$ in $D_{\bR^3}[0,\infty)$ as $n\to\infty$, where $W^1$ and $W^2$ are independent, standard one-dimensional Brownian motions.
\end{cor}

When $\{W_{a_n}\}$ and $\{W_{b_n}\}$ are such that $b_n/a_n \to L \in(0,\infty)$, the situation is much more delicate than in Theorem \ref{T:main2'}. We show that the function $\rho(t)$ is non zero, and in some cases it is not constant.

\begin{thm}\label{T:main2}
Let $\{a_n\}_{n=1}^\infty$ and $\{b_n\}_{n=1}^\infty$ be strictly increasing sequences in $\bN$. Let $L_n=b_n/a_n$ and suppose that $L_n\to L\in(0,\infty)$. Let $I=\{n:L_n=L\}$ and $c_n=\gcd(a_n,b_n)$.
  \begin{enumerate}[(i)]
  \item If $I^c$ is finite, then $L\in\bQ$ and, for all $t\ge0$,
    \[
    \lim_{n\to\infty}E[W_{a_n}(t)W_{b_n}(t)]
      = \frac{3t}{4p}\sum_{j=1}^q f_L(j/q),
    \]
  where $L=p/q$ and $p,q\in\bN$ are relatively prime.
  \item If $I$ is finite, then
    \[
    \lim_{n\to\infty}E[W_{a_n}(1)W_{b_n}(1)]
      = \frac3{4L}\int_0^1 f_L(x)\,dx.
    \]
  \item If $I$ is finite and $c_n\to\infty$, then
    \[
    \lim_{n\to\infty}E[W_{a_n}(t)W_{b_n}(t)]
      = \frac{3t}{4L}\int_0^1 f_L(x)\,dx,
    \]
  for all $t\ge0$.
  \item If there exists $k\in\bN$ such that $b_n=k\mod a_n$ for all $n$, then
    \[
    \lim_{n\to\infty}E[W_{a_n}(t)W_{b_n}(t)]
      = \frac3{4L}\int_0^t \wh f_L(kx)\,dx,
    \]
  for all $t\ge0$.
  \end{enumerate}
\end{thm}

\begin{rmk}
In Theorem \ref{T:main2} (iv), we assume that there exists $k\in\bN$ such that $b_n=k\mod a_n$ for all $n$. Note that this implies $k/c_n$ is an integer for all $n$. In particular, $\{c_n\}$ is a bounded sequence of integers. Moreover, since $\{c_n\}$ is bounded, this implies that $I$ is finite. Comparing Parts (ii) and (iv) of Theorem \ref{T:main2}, it follows that
  \[
  \frac3{4L}\int_0^1 \wh f_L(kx)\,dx
    = \frac3{4L}\int_0^1 f_L(x)\,dx,
  \]
for all $k\in\bN$. In fact, more can be said. Letting $y=kx$, we have
  \begin{align*}
  \frac3{4L}\int_0^t \wh f_L(kx)\,dx
    &= \frac3{4Lk}\int_0^{kt} \wh f_L(y)\,dy\\
  &= \frac3{4Lk}\bigg(\bigg(\sum_{j=1}^{\flr{kt}}
    \int_{j-1}^j \wh f_L(x)\,dx\bigg)
    + \int_{\flr{kt}}^{kt} \wh f_L(x)\,dx\bigg)\\
  &= \frac3{4Lk}\bigg(\flr{kt}\int_0^1 f_L(x)\,dx
    + \int_0^{kt-\flr{kt}} f_L(x)\,dx\bigg)\\
  &= \frac3{4L}\bigg(\frac{\flr{kt}}k\int_0^1 f_L(x)\,dx
    + \frac1k\int_0^{kt-\flr{kt}} f_L(x)\,dx\bigg).
  \end{align*}
Hence,
  \[
  \frac3{4L}\int_0^t \wh f_L(kx)\,dx
    = \frac{3t}{4L}\int_0^1 f_L(x)\,dx,
  \]
whenever $kt\in\bN$. \qed
\end{rmk}

 \noindent \textbf{Proof of Theorem \ref{T:main2}.} Let $\{a_n\}_{n=1}^\infty$ and $\{b_n\}_{n=1}^\infty$ be strictly increasing sequences in $\bN$. Let $L_n=b_n/a_n$ and suppose that $L_n\to L\in(0,\infty)$. Recall $\wt W_n(t)=W_n(t)-3n^{-1/3}B(\flr{nt}/n)$, and note that it will suffice to prove the corresponding limits for $\wt W$ rather than $W$.

Fix $t\in[0,1]$. Since $W_n(t)=0$ if $\flr{nt}=0$, we may assume $t>0$ and $n$ is sufficiently large so that $\flr{a_nt}>0$ and $\flr{b_nt}>0$. As in \eqref{main1.4}, we have
  \[
  S_n(t) := E[\wt W_{a_n}(t)\wt W_{b_n}(t)]
    = \frac34\sum_{j=1}^{\flr{a_nt}}\sum_{k=1}^{\flr{b_nt}}
    \Phi_n(j,k)^3.
  \]
Making the change of index $m=k-\flr{jL_n}$, we then have
  \[
  S_n(t) = \frac34\sum_{j=1}^{\flr{a_nt}}
    \sum_{m=1-\flr{jL_n}}^{\flr{b_nt}-\flr{jL_n}}
    \Phi_n(j,m + \flr{jL_n})^3.
  \]
Note that by \eqref{Phinfrel},
  \begin{equation}\label{main2.2}
  \Phi_n(j,m + \flr{jL_n})^3 
    = \frac1{b_n}f_{m,L_n}(jL_n - \flr{jL_n}).
  \end{equation}
If $1\le j\le\flr{a_nt}$, then
  \begin{align*}
  m \le \flr{b_nt} - \flr{jL_n}
    &\iff \flr{jL_n} < \flr{b_nt} - m + 1\\
  &\iff jL_n < \flr{b_nt} - m + 1\\
  &\iff j < \frac{\flr{b_nt} - m + 1}{L_n}\\
  &\iff j < \ceil{\frac{\flr{b_nt} - m + 1}{L_n}},
  \end{align*}
and also
  \begin{align*}
  m \ge 1 - \flr{jL_n} &\iff \flr{jL_n} \ge 1 - m\\
  &\iff jL_n \ge 1 - m\\
  &\iff j \ge \frac{1 - m}{L_n}\\
  &\iff j \ge \ceil{\frac{1 - m}{L_n}}.
  \end{align*}
Hence, when we reverse the order of summation, we obtain
  \[
  S_n(t) = \frac34\sum_{m=1-\flr{\flr{a_nt}L_n}}^{\flr{b_nt}-\flr{L_n}}
    \sum_{j=\ell_{m,n}}^{u_{m,n}}
    \Phi_n(j,m + \flr{jL_n})^3,
  \]
where
  \begin{align}
  \ell_{m,n} &= \ceil{\frac{1 - m}{L_n}} \vee 1,\label{ell_def}\\
  u_{m,n} &= \left({\ceil{\frac{\flr{b_nt} - m + 1}{L_n}} - 1}\right)
    \wedge \flr{a_nt}.\label{u_def}
  \end{align}
Let us define
  \[
  \be(m,n) = \sum_{j=\ell_{m,n}}^{u_{m,n}}
    \Phi_n(j,m + \flr{jL_n})^3,
  \]
so that we may write
  \[
  S_n(t) = \frac34\sum_{m\in\bZ}
    \be(m,n){\bf 1}_{[1 - \flr{\flr{a_nt}L_n},\flr{b_nt} - \flr{L_n}]}(m).
  \]
We wish to apply dominated convergence to this sum.

Choose an integer $M\ge2$ such that $L_n\le M$ for all $n\in\bN$. Define
  \[
  C_m = \begin{cases}
      8                  &\text{if $|m| \le M$},\\
      27(|m| - M)^{-3}   &\text{if $|m| > M$}.
    \end{cases}
  \]
We claim that $|\be(m,n)|\le C_mt$ for all $m$ and $n$. Once we prove this claim, we may use dominated convergence to conclude that
  \begin{equation}\label{main2.1}
  \lim_{n\to\infty} S_n(t) = \frac34\sum_{m\in\bZ}
    \lim_{n\to\infty}\be(m,n),
  \end{equation}
provided the limit on the right-hand side exists for each fixed $m$.

To prove the claim, first note that $1\le\ell_{m,n}\le u_{m,n}\le \flr{a_nt}$, so that
  \[
  |\be(m,n)| \le \sum_{j=1}^{\flr{a_nt}}
    |\Phi_n(j,m + \flr{jL_n})|^3.
  \]
Thus, by \eqref{eq2}, we have $|\be(m,n)| \le 8(a_n^{-1}\wedge b_n^{-1})\flr{a_nt} \le 8t$ for all $m$ and $n$. We therefore need only consider $|m|>M$.

First suppose $m>M$. Then
  \[
  \frac{m + \flr{jL_n} - 1}{b_n}
    > \frac{m + jL_n - 2}{b_n}
    > \frac{jL_n}{b_n} = \frac{j}{a_n}.
  \]
Hence, by \eqref{Phi}, \eqref{Phi_symmetry}, and \eqref{int_rep1},
  \begin{align*}
  |\Phi_n(j,m + \flr{jL_n})|
    &= \bigg|\frac29\int_0^{a_n^{-1}}\int_0^{b_n^{-1}}
    \left({\frac{m + \flr{jL_n} - 1}{b_n} - \frac j{a_n}
    + x + y}\right)^{-5/3}\,dx\,dy\bigg|\\
  &\le \int_0^{a_n^{-1}}\int_0^{b_n^{-1}}
    \left({\frac{m + jL_n - 2}{b_n} - \frac j{a_n}
    + y}\right)^{-5/3}\,dx\,dy\\
  &= \int_0^{a_n^{-1}}\int_0^{b_n^{-1}}
    \left({\frac{m - 2}{b_n} + y}\right)^{-5/3}\,dx\,dy.
  \end{align*}
By Lemma \ref{L:inequ},
  \[
  |\Phi_n(j,m + \flr{jL_n})|
    \le  \int_0^{a_n^{-1}}\int_0^{b_n^{-1}}
    \left({\frac{m - 2}{b_n}}\right)^{-1}y^{-2/3}\,dx\,dy
    = 3a_n^{-1/3}(m - 2)^{-1}.
  \]
Thus, $|\be(m,n)| \le 27a_n^{-1}\sum_{j=1}^{\flr{a_nt}}(m - 2)^{-3} 
\leq 27t(m - 2)^{-3} \le C_mt$.

Next, suppose $m<-M$. Then
  \[
  \frac{m + \flr{jL_n}}{b_n} \le \frac{m + jL_n}{b_n}
    < \frac{-L_n + jL_n}{b_n} = \frac{j - 1}{a_n}.
  \]
Hence, by \eqref{Phi} and \eqref{int_rep1},
  \begin{align*}
  |\Phi_n(j,m + \flr{jL_n})|
    &= \bigg|\frac29\int_0^{a_n^{-1}}\int_0^{b_n^{-1}}
    \left({\frac{j - 1}{a_n} - \frac{m + \flr{jL_n}}{b_n}
    + x + y}\right)^{-5/3}\,dx\,dy\bigg|\\
  &\le \int_0^{a_n^{-1}}\int_0^{b_n^{-1}}
    \left({\frac{j - 1}{a_n} - \frac{m + jL_n}{b_n}
    + y}\right)^{-5/3}\,dx\,dy\\
  &= \int_0^{a_n^{-1}}\int_0^{b_n^{-1}}
    \left({-\frac 1{a_n} - \frac m{b_n}
    + y}\right)^{-5/3}\,dx\,dy.
  \end{align*}
By Lemma \ref{L:inequ},
  \begin{multline*}
  |\Phi_n(j,m + \flr{jL_n})|
    \le  \int_0^{a_n^{-1}}\int_0^{b_n^{-1}}
    \left({-\frac 1{a_n} - \frac m{b_n}}\right)^{-1}
    y^{-2/3}\,dx\,dy\\
    = 3a_n^{-1/3}(-L_n - m)^{-1}
    = 3a_n^{-1/3}(|m| - L_n)^{-1}
    \le 3a_n^{-1/3}(|m| - M)^{-1}.
  \end{multline*}
Thus, $|\be(m,n)| \le 27a_n^{-1}\sum_{j=1}^{\flr{a_nt}}(|m| - M)^{-3} \le C_mt$. This proves our claim and establishes \eqref{main2.1}, provided the limit on the right-hand side exists for each fixed $m$.

Recalling \eqref{f_def}, let us now define
  \begin{equation}\label{main2.9}
  \wt\be(m,n) = \frac1{L_na_n}\sum_{j=1}^{\flr{a_nt}}
    f_{m,L}(jL_n - \flr{jL_n}).
  \end{equation}
We will first show that
  \begin{equation}\label{main2.3}
  \lim_{n\to\infty} |\wt\be(m,n) - \be(m,n)| = 0,
  \end{equation}
for each fixed $m\in\bZ$. Since $1\le\ell_{m,n}\le u_{m,n}\le\flr{a_nt}$, we have $(\wt\be-\be)(m,n) = A_{m,n} + B_{m,n}$, where
  \begin{align*}
  A_{m,n} &= \frac1{L_na_n}\sum_{j=1}^{\ell_{m,n}-1}
    f_{m,L}(jL_n - \flr{jL_n})
    + \frac1{L_na_n}\sum_{j=u_{m,n}+1}^{\flr{a_nt}}
    f_{m,L}(jL_n - \flr{jL_n}),\\
  B_{m,n} &= \sum_{j=\ell_{m,n}}^{u_{m,n}}
    \left({\frac1{L_na_n}f_{m,L}(jL_n - \flr{jL_n})
    - \Phi_n(j,m + \flr{jL_n})^3}\right).
  \end{align*}
By \eqref{f_basic}, we have
  \begin{equation}\label{main2.5}
  |f_{m,L}(x)| \le 8 \quad \text{for all $m$, $L$, and $x$}.
  \end{equation}
Thus,
  \[
  |A_{m,n}| \le \frac8{b_n}(\ell_{m,n} - 1 + \flr{a_nt} - u_{m,n}).
  \]
From \eqref{ell_def}, we see that $\limsup_{n\to\infty}\ell_{m,n}< \infty$. From \eqref{u_def}, we have
  \begin{align*}
  u_{m,n}
    &\ge \left({\frac{\flr{b_nt} - m + 1}{L_n} - 1}\right)
    \wedge \flr{a_nt}\\
  &= \flr{a_nt} + \left({\left({
    \frac{\flr{b_nt} - m + 1}{L_n} - 1 - \flr{a_nt}
    }\right)\wedge 0}\right)\\
  &= \flr{a_nt} - \left({\left({
    \flr{a_nt} - \frac{\flr{b_nt}}{L_n} + \frac{m - 1}{L_n} + 1
    }\right)\vee 0}\right).
  \end{align*}
Since
  \[
  \flr{a_nt} - \frac{\flr{b_nt}}{L_n}
    < a_nt - \frac{b_nt - 1}{L_n} = \frac1{L_n},
  \]
this gives
  \[
  \flr{a_nt} - u_{m,n} \le {\left({
    \frac m{L_n} + 1
    }\right)\vee 0},
  \]
and this shows that $\limsup_{n\to\infty}(\flr{a_nt}-u_{m,n})<\infty$. Hence, $A_{m,n}\to0$ as $n\to\infty$.

For $B_{m,n}$, we may use \eqref{main2.2} to write
  \[
  B_{m,n} = \frac1{b_n}\sum_{j=\ell_{m,n}}^{u_{m,n}}
    (f_{m,L} - f_{m,L_n})(jL_n - \flr{jL_n}).
  \]
By Lemma \ref{L:f_misc} (i),
  \[
  |f_{m,L}(x) - f_{m,L_n}(x)| \le 24|L_n - L|^{1/3},
  \]
for all $x$. This gives
  \[
  |B_{m,n}| \le \frac{24}{b_n}\sum_{j=1}^{\flr{a_nt}}|L_n - L|^{1/3}
    \le \frac{24t}{L_n}|L_n - L|^{1/3} \to 0,
  \]
as $n\to\infty$, and we have proved \eqref{main2.3}.

Finally, we calculate $\lim_{n\to\infty}\wt\be(m,n)$. We begin by rewriting $\wt\be(m,n)$ in the following way. For each $n$, choose $p_n,q_n,c_n\in\bN$ such that $a_n=c_nq_n$, $b_n=c_np_n$, and $p_n$ and $q_n$ are relatively prime. In general, if $p\in\bZ$ and $q\in\bN$, then let $[p]_q$ denote the unique integer such that $0\le[p]_q<q$ and $p\equiv[p]_q\mod q$. Note that
  \[
  [p]_q = q\left({\frac pq - \flr{\frac pq}}\right).
  \]
Thus,
  \begin{equation}\label{main2.8}
  jL_n - \flr{jL_n}
    = \frac{jb_n}{a_n} - \flr{\frac{jb_n}{a_n}}
    = \frac{jp_n}{q_n} - \flr{\frac{jp_n}{q_n}}
    = \frac1{q_n}[jp_n]_{q_n}.
  \end{equation}
Hence, by \eqref{main2.9},
  \begin{align*}
  \wt\be(m,n) &= \frac1{L_na_n}\sum_{j=1}^{\flr{a_nt}}
    f_{m,L}([jp_n]_{q_n}/q_n)\\
  &= \frac1{L_nc_nq_n}\sum_{j=1}^{\flr{c_nq_nt}}
    f_{m,L}([jp_n]_{q_n}/q_n).
  \end{align*}
Let $\al_n,r_n$ be the unique integers such that $\flr{c_nq_nt}=\al_n q_n+r_n$ and $0 \le r_n < q_n$. Note that $\al_n \ge 0$ and $r_n =[\flr{c_nq_nt}]_{q_n}$. Since $h\in\bZ$ implies $[p+hq]_q=[p]_q$, we have
  \begin{align*}
  \wt\be(m,n) &= \frac1{L_nc_nq_n}\bigg(
    \sum_{h=0}^{\al_n-1}\sum_{j=1}^{q_n}
    f_{m,L}([(j + hq_n)p_n]_{q_n}/q_n)
    + \sum_{j=1}^{r_n} f_{m,L}([(j + \al_n q_n)p_n]_{q_n}/q_n)
    \bigg)\\
  &= \frac{\al_n}{c_n}\frac1{L_nq_n}
    \sum_{j=1}^{q_n} f_{m,L}([jp_n]_{q_n}/q_n)
    + \frac1{L_na_n}\sum_{j=1}^{r_n} f_{m,L}([jp_n]_{q_n}/q_n).
  \end{align*}
Also note that if $p$ and $q$ are relatively prime, then
  \[
  \{[jp]_q: 1 \le j \le q\} = \{0,1,2,\ldots,q-1\}.
  \]
Therefore,
  \begin{equation}\label{main2.6}
  \wt\be(m,n) = \frac{\al_n}{c_n}
    \frac1{L_nq_n}\sum_{j=0}^{q_n-1} f_{m,L}(j/q_n)
    + \frac1{L_na_n}\sum_{j=1}^{r_n} f_{m,L}([jp_n]_{q_n}/q_n).
  \end{equation}

Now let $I=\{n:L_n=L\}$. First assume $I$ is finite and $t=1$. Then $r_n=0$ and $\al_n=c_n$, so that
  \[
  \wt\be(m,n) = \frac1{L_nq_n}\sum_{j=0}^{q_n-1} f_{m,L}(j/q_n).
  \]
We first prove that $\lim_{n\to\infty}q_n=\infty$. Let $M>0$ be arbitrary. Let $S =\{p/q: p\in\bZ,q\in\bN,q\le M\}$. Choose $\ep>0$ small enough so that $S\cap(L-\ep,L+\ep)\subset\{L\}$. Choose $n_0\in\bN$ large enough so that $I\subset\{1,\ldots,n_0\}$, and also $|L_n-L|<\ep$ for all $n>n_0$. Let $n>n_0$ be arbitrary. Then
  \[
  \frac{p_n}{q_n} = \frac{b_n}{a_n} = L_n
    \in (L - \ep, L + \ep) \setminus \{L\}.
  \]
Hence, $p_n/q_n\notin S$, which implies $q_n>M$, and this shows that $\lim_{n\to\infty}q_n=\infty$.

Since $f_{m,L}$ is continuous, it now follows that
  \begin{align*}
  \lim_{n\to\infty}\wt\be(m,n)
    &= \lim_{n\to\infty}
    \frac1{L_nq_n}\sum_{j=0}^{q_n-1} f_{m,L}(j/q_n)\\
  &= \frac1L\int_0^1 f_{m,L}(x)\,dx.
  \end{align*}
By \eqref{main2.1} and \eqref{main2.3}, we therefore have
  \[
  \lim_{n\to\infty} S_n(1) = \frac34\sum_{m\in\bZ}
    \frac1L\int_0^1 f_{m,L}(x)\,dx.
  \]
By Lemma \ref{L:fkLconv}, we may interchange  the summation and integration, and this proves Part (ii) of the theorem.

Next, assume $I^c$ is finite and $t>0$. In this case, there exists $n_0$ such that $L_n=L$ for all $n\ge n_0$. In particular, $L\in\bQ$, so we may write $L=p/q$, where $p,q\in\bN$ are relatively prime. In this case, $q_n=q$ for all $n\ge n_0$. Therefore, by \eqref{main2.6}, for all $n\ge n_0$,
  \[
  \wt\be(m,n) = \frac{\al_n}{c_n}
    \frac1{L_nq}\sum_{j=0}^{q-1} f_{m,L}(j/q) + \ep_n,
  \]
where, by \eqref{main2.5},
  \begin{equation}\label{main2.7}
  |\ep_n|
    = \bigg|\frac1{L_na_n}\sum_{j=1}^{r_n} f_{m,L}([jp]_q/q)\bigg|
    \le \frac{8r_n}{L_na_n}
    < \frac{8q}{L_na_n} \to 0,
  \end{equation}
as $n\to\infty$. Also,
  \[
  \left|{\frac{\al_n}{c_n} - t}\right|
    = \left|{\frac{\flr{c_nqt} - c_nqt - {r_n}}{c_nq}}\right|
    \le \frac{1 + q}{a_n} \to 0,
  \]
as $n\to\infty$. Thus,
  \[
  \lim_{n\to\infty}\wt\be(m,n)
    = \frac t{Lq}\sum_{j=0}^{q-1} f_{m,L}(j/q).
  \]
As above, using \eqref{main2.1} and \eqref{main2.3}, we have
  \[
  \lim_{n\to\infty}S_n(t)
    = \frac{3t}{4Lq}\sum_{m\in\bZ}\sum_{j=0}^{q-1} f_{m,L}(j/q)
    = \frac{3t}{4p}\sum_{j=0}^{q-1}\sum_{m\in\bZ} f_{m,L}(j/q).
  \]
Since $f_{m,L}(1)=f_{m-1,L}(0)$, we may write
  \[
  \lim_{n\to\infty}S_n(t)
    = \frac{3t}{4p}\sum_{j=1}^q\sum_{m\in\bZ} f_{m,L}(j/q)
    = \frac{3t}{4p}\sum_{j=1}^q f_L(j/q),
  \]
and this proves Part (i) of the theorem.

Next, assume $I$ is finite and $c_n\to\infty$. Note that
  \[
  r_n = [\flr{a_nt}]_{q_n} = q_n\left({
    \frac{\flr{a_nt}}{q_n} - \flr{\frac{\flr{a_nt}}{q_n}}
    }\right).
  \]
Thus,
  \[
  \al_n = \frac{\flr{a_nt} - r_n}{q_n}
    = \flr{\frac{\flr{a_nt}}{q_n}}.
  \]
It follows that $\al_n\le a_nt/q_n = c_nt$, and also
  \[
  \al_n > \frac{\flr{a_nt}}{q_n} - 1
    > \frac{a_nt - 1}{q_n} - 1
    = c_nt - \frac1{q_n} - 1.
  \]
Hence,
  \[
  t - \left({\frac1{a_n} + \frac1{c_n}}\right)
    < \frac{\al_n}{c_n} \le t.
  \]
Since both $c_n\to\infty$ and $a_n\to\infty$, this shows that $\al_n/c_n \to t$ as $n\to\infty$.

Also, as in \eqref{main2.7},
  \[
  \bigg|\frac1{L_na_n}\sum_{j=1}^{r_n} f_{m,L}([jp]_q/q)\bigg|
    < \frac{8q_n}{L_na_n}
    = \frac8{L_nc_n} \to 0,
  \]
as $n\to\infty$. Therefore, using \eqref{main2.6} and the argument immediately following \eqref{main2.6}, we have
  \[
  \lim_{n\to\infty}\wt\be(m,n) = \frac tL\int_0^1 f_{m,L}(x)\,dx.
  \]
By \eqref{main2.1} and \eqref{main2.3}, we therefore have
  \[
  \lim_{n\to\infty} S_n(t) = \frac{3t}4\sum_{m\in\bZ}
    \frac1L\int_0^1 f_{m,L}(x)\,dx.
  \]
By Lemma \ref{L:fkLconv}, we may interchange  the summation and integration, and this proves Part (iii) of the theorem.

Finally, assume there exists $k\in\bN$ such that $b_n=k\mod a_n$ for all $n$. As in \eqref{main2.8}, we may write $jL_n-\flr{jL_n} = [jb_n]_{a_n}/a_n$. Hence, by \eqref{main2.9},
  \[
  \wt\be(m,n) = \frac1{L_na_n}\sum_{j=1}^{\flr{a_nt}}
    f_{m,L}([jb_n]_{a_n}/a_n).
  \]
For $n$ sufficiently large, $k<a_n$, so that
  \[
  k = [b_n - a_n]_{a_n}
    = a_n\left({
    \frac{b_n-a_n}{a_n} - \flr{\frac{b_n-a_n}{a_n}}
    }\right).
  \]
Define $k_n=(b_n-k)/a_n$. Then $b_n = k_na_n+k$ and
  \[
  k_n = \frac{b_n}{a_n} - \left({
    \frac{b_n-a_n}{a_n} - \flr{\frac{b_n-a_n}{a_n}}
    }\right)
    = 1 + \flr{\frac{b_n-a_n}{a_n}} \in \bN.
  \]
Thus,
  \[
  [jb_n]_{a_n} = [j(k_na_n + k)]_{a_n} = [jk]_{a_n}
    = a_n\left({
    \frac{jk}{a_n} - \flr{\frac{jk}{a_n}}
    }\right),
  \]
giving
  \[
  \wt\be(m,n) = \frac1{L_na_n}\sum_{j=1}^{\flr{a_nt}}
    \wh f_{m,L}(jk/a_n).
  \]
Since $a_n\to\infty$ and $\wh f_{m,L}$ is continuous, we have
  \[
  \lim_{n\to\infty}\wt\be(m,n)
    = \frac1L\int_0^t \wh f_{m,L}(kx)\,dx.
  \]
By \eqref{main2.1} and \eqref{main2.3}, we therefore have
  \[
  \lim_{n\to\infty} S_n(t) = \frac34\sum_{m\in\bZ}
    \frac1L\int_0^t \wh f_{m,L}(kx)\,dx.
  \]
By Lemma \ref{L:fkLconv}, we may interchange  the summation and integration, and this proves Part (iv) of the theorem. \qed

\medskip

As a consequence of Theorem \ref{T:main2} we can establish the following result  on the convergence in distribution of the sequence $(B, W_{a_n} , W_{b_n})$.

\begin{cor}\label{C:main3}
Let $\{a_n\}_{n=1}^\infty$ and $\{b_n\}_{n=1}^\infty$ be strictly increasing sequences in $\bN$. Let $L_n=b_n/a_n$ and suppose that $L_n\to L\in(0,\infty)$. Let $I=\{n:L_n=L\}$ and $c_n=\gcd(a_n,b_n)$. Given $\rho\in C[0,\infty)$, let $\si$ be given by \eqref{main1a} and $X^\rho$ by \eqref{main1b}.
  \begin{enumerate}[(i)]
  \item If $I^c$ is finite, then $(B,W_{a_n},W_{b_n}) \To (B,X^\rho)$ in $D_{\bR^3}[0,\infty)$ as $n\to\infty$, where
    \begin{equation}\label{main2i}
    \rho(t) = \frac3{4p}\sum_{j=1}^q f_L(j/q),
    \end{equation}
  for all $t\ge0$. Here, $L\in\bQ$ and $p$ and $q$ are determined by $L=p/q$, where $p,q\in\bN$ are relatively prime.
  \item If $I$ is finite and $c_n\to\infty$, then $(B,W_{a_n},W_{b_n}) \to (B,X^\rho)$ in law in $D_{\bR^3}[0,\infty)$ as $n\to\infty$, where
    \[
    \rho(t) = \frac3{4L}\int_0^1 f_L(x)\,dx,
    \]
  for all $t\ge0$.
  \item If there exists $k\in\bN$ such that $b_n=k\mod a_n$ for all $n$, then $(B,W_{a_n},W_{b_n}) \to (B,X^\rho)$ in law in $D_{\bR^3}[0,\infty)$ as $n\to\infty$, where
    \[
    \rho(t) = \frac3{4L}\wh f_L(kt),
    \]
  for all $t\ge0$.
  \end{enumerate}
\end{cor}

\pf  First assume $I^c$ is finite. Let $0\le s\le t<\infty$. Let $\rho$ be given by \eqref{main2i}. By Theorem \ref{T:main2} and Remark \ref{R:main},
  \[
  \lim_{n\to\infty}E[\wt W_{a_n}(s)\wt W_{b_n}(s)]
    = \int_0^s \rho(u)\,du.
  \]
By Lemma \ref{L:main1}, this gives
  \[
  \lim_{n\to\infty}E[\wt W_{a_n}(s)\wt W_{b_n}(t)]
    = \int_0^s \rho(u)\,du.
  \]
Part (i) of the theorem now follows from Theorem \ref{T:main1} and Remark \ref{R:main}. The proofs of Parts (ii) and (iii) are similar. \qed

\section{Remarks and examples}

Let $\{a_n\}_{n=1}^\infty$ and $\{b_n\}_{n=1}^\infty$ be strictly increasing
sequences in $\bN$. For each $t>0$, let
  \[
  \ga^{a,b}(t) = \lim_{n\to\infty}
    \frac{E[W_{a_n}(t)W_{b_n}(t)]}{\ka^2 t},
  \]
provided this limit exists. Then $\ga^{a,b}(t)$ is the asymptotic correlation
of $W_{a_n}(t)$ and $W_{b_n}(t)$ as $n\to\infty$. Under the hypotheses of Corollary \ref{C:main3}, we have
  \begin{equation}\label{ga_rho}
  \ga^{a,b}(t) = \frac1{\ka^2t}\int_0^t \rho(x)\,dx.
  \end{equation}
Note that $\ga^{a,b}$ is a constant function if and only $\rho$ is constant, as in Corollary \ref{C:main3} (i) and (ii). Also note that, since $\wh f_L(k\,\cdot)$ is not a constant function by Lemma \ref{L:f_misc}, Corollary \ref{C:main3} (iii) shows that there are circumstances under which $\ga^{a,b}$ is not constant.

\begin{expl} \label{ex4.1}
If $a_n=b_n=n$ for all $n$, then $E[W_{a_n}(t)W_{b_n}(t)]\to\ka^2t$ as $n\to \infty$, so that $\ga^{a,b}\equiv1$. By \eqref{kap_def} and \eqref{f_def}, we observe that
  \begin{equation}\label{ka_f1}
  \ka^2 = \frac34\sum_{m\in\bZ} f_{m,1}(0) = \frac34 f_1(0).
  \end{equation}
Equivalently, we may write
  \[
  \ka^2 = 6\sum_{m\in\bZ}(E[(B(1) - B(0))(B(m + 1) - B(m))])^3.
  \]
For $L\in\bN$, let us define
  \[
  \ka_L^2 = \frac6L\sum_{m\in\bZ}
    (E[(B(1) - B(0))(B(m + L) - B(m))])^3.
  \]
Then using \eqref{Phi_def}, \eqref{f_def}, and Lemma \ref{L:f_misc} (ii), we have
  \[
  \ka_L^2 = \frac3{4L}\sum_{m\in\bZ}
    \Phi(0,1,m,m+L)^3 = \frac3{4L}f_L(0) = \frac3{4L}f_L(1).
  \]
If $a_n=n$ and $b_n=Ln$, then using Theorem \ref{T:main2} (i) with $p=L $ and $q=1$, as well as Lemma \ref{L:f_misc} (ii) and \eqref{f_def}, we obtain $\lim_{n\to\infty} E[W_n(t) W_{Ln}(t)] = \ka_L^2t$, giving $\ga^{a,b}(t) \equiv\ka_L^2/\ka^2$.

Numerical calculations suggest that, in this family of examples, $\ga^{a,b}$ decreases fairly quickly with $L$. For example, when $L=2$, we have $\ga^{a,b} \approx0.201928$, and when $L=5$, we have $\ga^{a,b}\approx0.043837$. \qed
\end{expl}

The scaling property of fBm manifests itself in the present investigation via the following result.

\begin{lemma}\label{L:scaling}
Let $\{a_n\}_{n=1}^\infty$ and $\{b_n\}_{n=1}^\infty$ be strictly increasing sequences in $\bN$, and $r\in(0,\infty)$. Assume $a_n^*= ra_n$ and $b_n^*=rb_n$ are integers for all $n$. Fix $t>0$. Then
  \[
  \lim_{n\to\infty} E[W_{a_n}(rt)W_{b_n}(rt)]
    = r\lim_{n\to\infty} E[W_{a_n^*}(t)W_{b_n^*}(t)],
  \]
provided that one of the two limits exist.
\end{lemma}

\pf Recall $\wt W_n(t)=W_n(t)-3n^{-1/3}B(\flr{nt}/n)$, and note that it will suffice to prove the lemma for $\wt W$ rather than $W$. As in \eqref{main1.4}, we have
  \[
  E[\wt W_{a_n}(rt)\wt W_{b_n}(rt)]
    = \frac34\sum_{j=1}^{\flr{a_nrt}}\sum_{k=1}^{\flr{b_nrt}}
    \Phi_n^{a,b}(j,k)^3.
  \]
Note that $\Phi_n^{a^*,b^*}=r^{-1/3}\Phi_n^{a,b}$. Thus,
  \begin{multline*}
  E[\wt W_{a_n}(rt)\wt W_{b_n}(rt)]
    = \frac{3r}4\sum_{j=1}^{\flr{a_nrt}}\sum_{k=1}^{\flr{b_nrt}}
    \Phi_n^{a^*,b^*}(j,k)^3\\
    = \frac{3r}4\sum_{j=1}^{\flr{a_n^*t}}\sum_{k=1}^{\flr{b_n^*t}}
    \Phi_n^{a^*,b^*}(j,k)^3
    = rE[\wt W_{a_n^*}(t)\wt W_{b_n^*}(t)].
  \end{multline*}
Letting $n\to\infty$ completes the proof. \qed

\begin{expl}
At first glance, Lemma \ref{L:scaling} may seem to suggest that $\ga^{a,b}$ should always be the constant function $\ga^{a,b}(t)\equiv\ka^{-2}E[W_{a_n}(1) W_{b_n}(1)]$. But, of course, we know this to be false from Corollary \ref{C:main3} (iii). A simple example illustrating this is the following.

Fix $L,k\in\bN$. Let $a_n=n$ and $b_n=Ln+k$. Note that
  \[
  L_n = \frac{b_n}{a_n} = L + \frac kn \to L,
  \]
as $n\to\infty$. By Corollary \ref{C:main3} (iii), \eqref{ga_rho}, and \eqref{ka_f1},
  \begin{equation}\label{ga_n+k}
  \ga^{a,b}(t) = \frac1{Lf_1(0)t}\int_0^t \wh f_L(kx)\,dx.
  \end{equation}
Lemma \ref{L:f_misc} shows that this is not a constant function, at least when $L=1$. In this example, Lemma \ref{L:scaling} implies that for $r\in\bN$,
  \[
  \lim_{n\to\infty} E[W_n(rt)W_{Ln+k}(rt)]
    = r\lim_{n\to\infty} E[W_{rn}(t)W_{rLn+rk}(t)].
  \]
This does not contradict \eqref{ga_n+k}, since $\{(W_{rn},W_{rLn+rk})\}_{n=1}^\infty$ is not a subsequence of $\{(W_n,W_{Ln+k})\}_{n=1}^\infty$. Rather, it is a subsequence of $\{(W_n,W_{Ln+rk})\}_{n=1}^\infty$. Hence, Lemma \ref{L:scaling}
is, in this case, illustrating the equality,
  \[
  \int_0^{rt} \wh f_L(kx)\,dx
    = r\int_0^t \wh f_L(rkx)\,dx,
  \]
which is easily verified by a simple change of variable.

Another interesting feature illustrated here is the following. If we fix $L=1$, then we have a family of examples indexed by $k$ that all share the same limiting ratio, $L$, yet produce different asymptotic correlation functions. Indeed, if it were the case that $\wh f_1(k_1\,\cdot) = \wh f_1(k_2\, \cdot)$ for some $k_1<k_2$, then we would have $\wh f_1(1/2)=\wh f_1(k_1^n/2k_2^n)$ for all $n$, and by continuity, $\wh f_1(1/2)=\wh f_1(0)$, contradicting Lemma \ref{L:f_misc}.

Note that, by the continuity of $f_L$, we have $\ga^{a,b}(t)\to f_L(0)/(Lf_1(0))$ as $t\downarrow 0$. Numerical calculations for the case $L=k=1$ suggest that $\ga_{a,b}$ is a positive function with $\ga^{a,b}(0.8) \approx 0.0750475$, so that the asymptotic correlation between $W_n(t)$ and $W_{n+1}(t)$ varies dramatically with $t$.
\qed
\end{expl}

\begin{expl} \label{ex4.2}
As an example illustrating Corollary \ref{C:main3} (ii), let $L\in\bN$ and consider $a_n=n^2$ and $b_n=Ln^2+n$. Then $L_n=b_n/a_n=L+1/n\to L$. Since $c_n =\gcd(a_n,b_n)= n$, Corollary \ref{C:main3} (ii), \eqref{ga_rho}, and \eqref{ka_f1} give
  \[
  \ga^{a,b}(t) \equiv \frac1{Lf_1(0)}\int_0^1 f_L(x)\,dx.
  \]
Numerical calculations suggest that for $L=1$ and $L=2$, $\ga^{a,b} \approx 0.101932$ and $\ga^{a,b} \approx 0.0468229$, respectively. Note that these numbers are several times smaller than the corresponding numbers for the sequences $a_n=n^2$ and $b_n=Ln^2$, which are covered by Example \ref{ex4.1}.
\qed
\end{expl}

\begin{expl}
Our penultimate example illustrates a situation where $c_n=\gcd(a_n,b_n)$ is constant, yet the asymptotic correlation $\ga^{a,b}$ does not exist.

Fix $k\in\bN$. Let $a_n=kn^2$ and
  \[
  b_n = \begin{cases}
    kn^2 + 2k &\text{if $n$ is odd},\\
    kn^2 + k &\text{if $n$ is even}.
    \end{cases}
  \]
Then $c_n=\gcd(a_n,b_n)=k$ for all $n$. By Theorem \ref{T:main2} (iv),
  \[
  \lim_{\substack{n\to\infty\\n\text{ odd}}}
    E[W_{a_n}(t)W_{b_n}(t)] = \lim_{m\to\infty}
    E[W_{k(2m+1)^2}(t)W_{k(2m+1)^2+2k}(t)]
    = \frac34\int_0^t \wh f_1(2kx)\,dx.
  \]
On the other hand,
  \[
  \lim_{\substack{n\to\infty\\n\text{ even}}}
    E[W_{a_n}(t)W_{b_n}(t)] = \lim_{m\to\infty}
    E[W_{4km^2}(t)W_{4km^2+k}(t)]
    = \frac34\int_0^t \wh f_1(kx)\,dx.
  \]
Since these are different functions, the sequence $(W_{a_n},W_{b_n})$ does
not converge and $\ga^{a,b}$ does not exist.
\qed
\end{expl}

\begin{expl}
Finally, we collect what might be called some non-examples, a few cases which are not covered by our present results. The first is $a_n=n^2$ and $b_n=(n+1)^2$. In this case, $L_n\to 1$, but $\gcd(a_n,b_n)=1$ for all $n$, and $b_n-a_n=2n+1 \mod a_n$. By Theorem \ref{T:main2} (ii), we know that $\ga^{a,b}(1)$ exists, but the existence and value of $\ga^{a,b}(t)$ for $t\ne1$ is not covered by our results.

The second non-example is $a_n=2n$ and $b_n=3n+1$. In this case, $L_n\to 3/2$, but $\gcd(a_n,b_n) \le 2$ for all $n$, and $b_n-a_n=n+1\mod a_n$. Again our results fail to give a complete picture of the function $\ga^{a,b}$.

Our last non-example is the following. Let $\al\in(1,2)$ be an irrational number whose decimal expansion contains only the digits 1, 3, 7, and 9. In other words, $\al=1 + \sum_{j=1}^\infty c_j10^{-j}$, where $c_j\in\{1,3,7,9\}$ for all $j$. Let $s_n=\sum_{j=1}^n c_j10^{-j}$, and define $a_n=10^n$ and $b_n=10^n(1+s_n)$. In this case, $L_n\to\al$, but $\gcd(a_n,b_n)=1$ for all $n$, and $b_n-a_n=10^ns_n\mod a_n$. Again our results tell us only the existence and value of $\ga^{a,b}(1)$.

There are, of course, many examples such as these which are not covered by Theorem \ref{T:main2}. Developing a more general set of results that describe the asymptotic behavior of the correlation of $W_{a_n}(t)$ and $W_{b_n}(t)$ in these examples is an open problem to be studied in subsequent work. \qed
\end{expl}

\section{Appendix}
In this section we include  a couple of  technical results that are used along the paper.

\begin{lemma}\label{L:inequ}
If $a_j$ and $p_j$ are positive real numbers, then
  \[
  \bigg(\sum_{j=1}^n a_j\bigg)^{-\sum_{j=1}^n p_j}
    \le \prod_{j=1}^n a_j^{-p_j}.
  \]
\end{lemma}

\pf  For every $k=1,\dots, n$ we have
\[
\left(\sum_{j=1}^n a_j \right) ^{-p_k} \le a_k^{-p_k},
\]
and the desired result follows by taking the product of these terms in $k$. \qed

\begin{lemma}\label{L:numeric}
Let $a,b\in\bZ$ with $a\le b$. Let $C \in \bR$ and $p>1$. If $b+1<C$, then
  \[
  \sum_{j=a}^b (C - j)^{-p} \le \frac1{p-1}(C-(b+1))^{-(p-1)}.
  \]
If $C<a-1$, then
  \[
  \sum_{k=a}^b (k - C)^{-p} \le \frac1{p-1}(a - 1 - C)^{-(p-1)}.
  \]
\end{lemma}

\pf If $b+1<C$, then
  \begin{multline*}
  \frac1{p-1}(C-(b+1))^{-(p-1)}
    \ge \frac1{p-1}((C-(b+1))^{-(p-1)} - (C-a)^{-(p-1)})\\
  = \int_a^{b+1} (C - x)^{-p}\,dx
    = \sum_{j=a}^b \int_j^{j+1} (C - x)^{-p}\,dx
    \ge \sum_{j=a}^b (C - j)^{-p}.
  \end{multline*}
If $C<a-1$, then
  \begin{multline*}
  \frac1{p-1}(a - 1 - C)^{-(p-1)}
    \ge \frac1{p-1}((a-1-C)^{-(p-1)} - (b-C)^{-(p-1)})\\
    = \int_{a-1}^b (x - C)^{-p}\,dx
    = \sum_{k=a}^b \int_{k-1}^k (x - C)^{-p}\,dx
    \ge \sum_{k=a}^b (k - C)^{-p}.
  \end{multline*}
\qed



\end{document}